\documentclass[11pt, reqno]{amsart}
\usepackage[utf8]{inputenc}
\usepackage{amsthm}
\usepackage{amsmath,amsfonts,amssymb}
\usepackage[a4paper, margin=1in]{geometry}
\usepackage{bm}
\usepackage{siunitx}
\usepackage{cite}
\usepackage{mwe}
\usepackage{graphicx}
\usepackage{caption}
\usepackage{subcaption}
\usepackage{hyperref}
\usepackage{soul}

\newtheorem{remark}{Remark}[section]

\newtheorem{theorem}{Theorem}[section]

\newcommand{\bld}{\mathbf}
\def\n{{n}}
\def\nm{{n\num{-1}}}

\newcommand{\quotes}[1]{``#1''}

\begin{document}
\title[FEM for PME]{Two Finite Element Approaches For The Porous Medium Equation That Are Positivity Preserving And Energy Stable}
\author{Arjun Vijaywargiya}
\author{Guosheng Fu}
\address{Department of Applied and Computational Mathematics and
Statistics, University of Notre Dame, USA.}
\email{avijaywa@nd.edu, gfu@nd.edu}
 \thanks{This research is partially supported by NSF grant DMS-2012031. }

\begin{abstract}
  In this work, we present the construction of two distinct finite element approaches to solve the Porous Medium Equation (PME). In the first approach, we transform the PME to a log-density variable formulation and construct a continuous Galerkin method. In the second approach, we introduce additional potential and velocity variables to rewrite the PME into a system of equations, for which we construct a mixed finite element method. Both approaches are first-order accurate, mass conserving,  and proved to be unconditionally energy stable for their respective energies. The mixed approach is shown to preserve positivity under a CFL condition, while a much stronger property of unconditional bound preservation is proved for the log-density approach. A novel feature of our schemes is that they can handle compactly supported initial data without the need for any perturbation techniques. Furthermore, the log-density method can handle unstructured grids in any number of dimensions, while the mixed method can handle unstructured grids in two dimensions. We present results from several numerical experiments to demonstrate these properties. 
\end{abstract}\maketitle

\section{Introduction}
The porous medium equation (PME) is given by 
\begin{equation}
    \label{PME}
    \frac{\partial \rho}{\partial t}-\nabla \cdot \nabla\left(\rho^m\right)=0,
\end{equation}
where $x \in \Omega \subset \mathbb R^d$, $\rho(t,x) \geq 0$ is the unknown density function, and $m \geq 1$. The PME can be found in many physical and biological phenomena, including the flow of an ideal gas through porous media, groundwater infiltration, the spread of viscous fluids, and boundary layer theory \cite{Gratton98,
Vazquez07, b03/04}. 
The PME is an example of a non-linear parabolic PDE with an underlying Wasserstein gradient flow structure \cite{Otto01} that guarantees that any solution for initially positive data will remain so for all times, and will be energy dissipative. 
Moreover, the PME has the peculiar \textit{finite speed of propagation} property which states that any compactly supported initial data will remain so for all times, and thus, the solution at the domain interface, referred to as the free boundary, will only propagate at a finite speed. The PME admits a classical self-similar weak solution \cite{Barenblatt1952,10.1093/qjmam/12.4.407} that is compactly supported and shows this property. It is also a commonly known fact that for certain initial data, the solution exhibits the \textit{waiting-time phenomenon} in which the free-boundary does not move, but the interior profile continues to evolve until a certain positive finite time \cite{doi:10.1137/0514049}.

The existence of a free boundary due to the degeneracy of the PME makes traditional parabolic numerical techniques ineffective. Various numerical methods have been developed for the PME \cite{Graveleau71,Jin98,Zhang09,Westdickenberg10, Chen16, Ngo17, Duan19, Liu20a,GU2020109378,yuanliu11,Duque13,Pop2002,Monsaingeon16,Gurtin1984,Bertsch1990,Benedetto1984}. We summarize some of the broad themes of these methods:
\begin{itemize}
    \item Standard numerical routines that invoke standard finite element procedures for spatial discretization and a predictor-corrector formulation for temporal discretization (PCSFE method) can suffer from oscillations at the free boundary which cannot be suppressed by increasing the polynomial degree or refining the spatial mesh \cite{Zhang09}.

    \item Discontinuous Galerkin \cite{Chen16}, Local Discontinuous Galerkin \cite{Zhang09}, and WENO \cite{yuanliu11} methods have been adapted to the PME to suppress these oscillations arising near the free boundary. Particularly, \cite{Zhang09} makes use of a non-negativity preserving limiter, while \cite{Chen16} implements a maximum-principle-satisfying (MPS) limiter. However, these schemes introduce additional numerical viscosity, due to which they fail to accurately track the free boundary, and thus, to accurately estimate the waiting time.

    \item Several interface tracking schemes \cite{Monsaingeon16,Gurtin1984,Bertsch1990,Benedetto1984} have been developed that track the free boundary by solving the equation of the interface in the Lagrangian coordinate. However, these schemes have limited applicability in higher dimensions and to initial data with complex support due to the complexity of implementation. Recent works in \cite{Carrillo,liu20} develop fully Lagrangian schemes 
    that can be applied to 2D problems. Nonetheless, being a Lagrangian scheme, it may not handle well initial data with complex support. If any \quotes{tangling} of the mesh occurs, the solution has to be manually interpolated onto a new mesh.

    \item Perturbation techniques have also been used to remove the degeneracy of the PME \cite{Duque13,Pop2002,GU2020109378}. In particular, the non-negative initial data is perturbed by a small parameter $\varepsilon > 0$ to make it positive everywhere. Using this approach, \cite{GU2020109378} constructs a first-order scheme that is unconditionally energy dissipative as well as provably bound preserving on structured tensor-product grids. However, in addition to the perturbation by $\varepsilon$, this scheme has the limitation that it cannot be applied to complex geometries.
\end{itemize}

In this work, we present two positive and energy-preserving nonlinear schemes for the PME: one is based on a log-density formulation, and the other uses a mixed formulation. Our log-density formulation scheme is closely related to the work in \cite{GU2020109378} but can handle fully unstructured grids in any number of dimensions. Additionally, both schemes do not require any perturbation of a non-positive initial data. 
\begin{itemize}
    \item The log-density based scheme uses a classical conforming piecewise linear finite element space for the spatial discretization in combination with a first-order semi-implicit time discretization. 
Global mass conservation, the positivity of density, and energy dissipation for the energy $\int_\Omega\rho(\log(\rho)-1)\,\mathrm{dx}$ are proven for this fully discrete scheme on general unstructured meshes. Unique solvability of the nonlinear equations in each time step is established when the mass matrix is lumped to be a diagonal matrix. 
Moreover, under the condition of a Delaunay triangulation along with a special edge-based discretization of the nonlinear diffusion term \cite{Xu1999}, the scheme can be further proven to satisfy a discrete maximum principle.
While the theoretical results for this scheme are proven under the condition that the initial density $\rho^0(x)>0$ is positive everywhere, we show a practical implementation to handle compactly supported
initial data which simply deactivates the degrees of freedom if the associated diagonal entry of the stiffness matrix is below a threshold value (e.g., $10^{-14}$).
\item 
On the other hand, the mixed method is based on a first-order reformulation of the PME, where the unknowns are density, potential, and velocity. 
The classical RT0-P0 finite element pair is used for spatial discretization, where 
the velocity is approximated via the lowest order Raviart-Thomas finite element space \cite{RT}, and the
density and potential are approximated by a discontinuous piecewise constant space. 
A first-order semi-implicit time discretization is then applied to this spatial discretization, which results in a nonlinear fully discrete scheme.
The fully discrete scheme is locally mass conservative and positivity preserving under a classical CFL condition on general unstructured meshes. 
Furthermore, we prove energy dissipation 
for the {\it physical energy} $\int_{\Omega}\frac{\rho^m}{m-1}\,\mathrm{dx}$
of this scheme 
when mass lumping \cite{Baranger1996} is applied to evaluate the velocity mass matrix. 
Positivity of this lumped velocity mass matrix requires the triangulated mesh to be {Delaunay}.
We note that this mass lumping is simply a trapezoidal rule on tensor-product meshes, on which one can further prove the solvability of the nonlinear system and unconditional positivity preservation following similar arguments as in the log-density formulation; see also \cite{GU2020109378}. 
The nonlinear system in each time step can be efficiently solved by expressing velocity and potential in terms of the density unknowns
and then solving the parabolic system of nonlinear equations for density alone 
using Newton's method.
This scheme naturally handles a compactly supported initial density profile as zero density/potential is allowed in the scheme.
\end{itemize}



The rest of the paper is organized as follows. In section 2, we describe the construction of the log-density scheme and prove several properties like mass conservation, energy stability, unique solvability, and bound preservation. In section 3, we describe the mixed method and prove properties like local mass dissipation, energy stability, and positivity. In section 4, we discuss several numerical experiments and compare the results of the two schemes. In section 5, we provide some concluding remarks. 

\section{Log-Density Formulation}
We consider the PME \eqref{PME} expressed using the log-density variable $u := \log(\rho)$ on a bounded polyhedral domain $\Omega \subset \mathbb R^d, d=1,2,3$ with a homogeneous Neumann boundary condition:
\begin{subequations}
   \label{pmeX}
\begin{align}
  \label{pme1}
  \frac{\partial \exp(u)}{\partial t}-\nabla\cdot
  \left( m\exp(m\,u) \nabla u\right)&=0,\quad {\text{in} }\quad \Omega,
  \\
  \label{pme2}
  \frac{\partial u}{\partial n} &= 0, \quad {\text on }\quad \partial\Omega,
\end{align}
and initial data
 \begin{alignat}{2}
  \label{pme3}
  u(0,x) =&\; \log(\rho^0(x)) \text{ in }\Omega,
\end{alignat}
 \end{subequations}
 where the initial density 
 $\rho^0(x)>0$ is assumed to be positive everywhere in the domain. 
 This {\it log-density} based formulation was first developed in \cite{Metti16}
for the Poisson-Nernst-Planck equations, see also \cite{FuXu21a}, 
which is closely related to the 
the entropy-stable schemes based on  the {\it entropy
variables} for hyperbolic  conservation laws and compressible flow in the CFD literature
\cite{Harten83, Tadmor84,Hughes86_2,Barth99}.
In particular, the density $\rho = \exp(u)$ is guaranteed to stay positive as long as the initial data 
$\rho^0(x)$ is positive.
We note that the recently introduced bound preserving and energy dissipative finite difference scheme by Gu and Shen \cite{GU2020109378} is also closely related to this log-density formulation.

Note that the PME has a special \textit{finite speed of propagation} property \cite{Vazquez07} which states that if the initial data $\rho^0$ has compact support, then the solution to the Cauchy problem of the PME will also have compact support at any other time, $t > 0$. For this reason, special care has to be taken in the case of compactly supported initial data as $u$ is negative infinity at locations where $\rho^0$ is zero. In \cite{GU2020109378}, the authors add a small perturbation $\varepsilon$ to the initial data to make it positive everywhere. We employ a different approach in which we deactivate degrees of freedom wherever  $\rho^0$ is close to zero; see more discussion at the end of this section.

The PME \eqref{pmeX} satisfy the following three important properties:
\begin{subequations}
  \label{prop} 
\begin{itemize}
  \item [(i)] Mass conservation: 
    \begin{align}
      \label{mass}
      \int_{\Omega}\rho(t, x)dx = 
    \int_{\Omega}\rho^0(x)dx. 
    \end{align}
  \item [(ii)] Positivity: 
\begin{align}
  \label{positivity}
  \text{
  If $\rho^0(x)>0$, then $\rho(t, x)>0$ for any $t>0$.}
\end{align}
  \item [(iii)] Energy dissipation:
    \begin{align}
      \label{ener}
      \frac{d}{dt} E=-\int_{\Omega}m\rho^m|\nabla u|^2dx,
    \end{align}
    where the energy  $E$ is given by
\begin{align}
  \label{enerXX}
E:= \int_{\Omega}\rho(\log(\rho)-1)dx
= \int_{\Omega}\exp(u)(u-1)dx.
\end{align}   
\end{itemize}
\end{subequations}
Our goal is to design a numerical scheme that preserves 
these three properties.

\subsection{Spatial and Temporal Discretizations}
\label{sec:notation}
We describe the method on a general unstructured simplicial mesh, although quadrilateral/hexahedral meshes can also be used.
Let $\Omega_h:=\{K_i\}_{i=1}^{N_K}$ be a conforming simplicial
triangulation of the domain $\Omega$ with $N_K$ elements.
Denote $\mathcal{E}_h=\{E_i\}_{i=1}^{N_E}$ as the collection of $N_E$ edges of 
$\Omega_h$, and $\mathcal{V}_h=\{v_i\}_{i=1}^{N_V}$ as its collection 
of $N_V$ vertices.
For any element $K\in \Omega_h$, denote
$\mathcal{E}_K:=\{E\in \mathcal{E}_h: \; E\subset \bar K\}$ as its edges, and 
$\mathcal{V}_K:=\{v\in \mathcal{V}_h: \; v\in \bar K\}$ as its vertices.

We shall use the $H^1$-conforming finite element space
 \begin{align}
   \label{space-pk}
   V_h := \{v_h\in H^1(\Omega):\; v_h|_K\in P_1 (K),\quad \forall K\in\Omega_h\},
 \end{align}
 where $P_1(K)$ is the space of linear polynomials on a simplex $K$. 
The space $V_h$ is equipped with the standard nodal (hat) basis $\{\phi_{i}(x)\}_{i=1}^{N_V}$ in which $\phi_{i}(v_j) = \delta_{ij}$ where 
$\delta_{ij}$ is the Kronecker delta function.
Hence any function $w_h\in V_h$ can be expressed as 
\[
w_h=\sum_{i=1}^{N_V}w_i\phi_i, 
\]
where $\underline{w}:=[w_1,\cdots, w_{N_V}]'$
is the coefficient vector satisfying
$w_i=w_h(v_i)$.

 The spatial discretization for \eqref{pmeX} then reads: find $u_h \in V_h$ such that, for $t>0$,
 \begin{align}
    \label{fem}
    \int_{\Omega}\frac{\partial \exp(u_{h})}{\partial t} v_h\,\mathrm{dx} 
    +\int_{\Omega} 
    m\exp(m\,u_{h})\nabla u_{h}\cdot\nabla v_h\,\mathrm{dx}
=&0,\quad 
\forall v_h\in V_h,
  \end{align}
 with initial conditions
 \begin{equation}
 \label{ic}
u_h(0, v_i)=\log \left(\rho^0_h(v_i)\right), \quad \forall v_i\in \mathcal{V}_h.
\end{equation}

We apply a first-order semi-implicit discretization for the ODE system \eqref{fem} to arrive at the following fully discrete scheme:
Given data $u_h^{n-1}\in V_h$ at time $t^{n-1}$
and time step size $\Delta t$,
find $u_h^\n \in V_h$ at time $t^n=t^{n-1}+\Delta t$ such that
\begin{equation}
\label{semi-imp}
\mathcal{M}\left(\frac{\exp (u_{h}^\n)-\exp (u_h^\nm)}{\Delta t}, v_h\right)+
\mathcal{A}\left(m\exp (m u_h^\nm); 
\nabla u_h^\n, \nabla v_h\right)=0, \quad \forall v_h \in V_h,
\end{equation}
where the mass operator $\mathcal{M}$ and stiffness operator $\mathcal{A}$ read as follows:
\begin{align*}
\mathcal{M}(\alpha, \beta):=&\; 
\int_{\Omega} \alpha\cdot\beta\mathrm{dx},\quad\quad
\mathcal{A}(\gamma; \alpha, \beta):=\; 
\int_{\Omega} \gamma\nabla\alpha\cdot\nabla \beta\mathrm{dx}.
\end{align*}


Unique solvability of the scheme requires the mass matrix to be mass lumped, which is achieved by applying the following vertex-based quadrature rule
\begin{equation}
    \label{mlquad}
\mathcal{M}_h(\alpha, \beta) := 
\sum_{K\in\Omega_h} \sum_{v\in \mathcal{V}_K}\frac{|K|}{d+1} \alpha(v) \beta(v)=
\sum_{i=1}^{N_v}\frac{|S_i|}{d+1} \alpha(v_i) \beta(v_i),
\end{equation}
where $S_i:=\cup_{v_i\in \bar K}\{\bar K\}$ is the vertex patch of $v_i$ and 
$|S_i|$ is its volume. 

Furthermore, we make use of the following edge-based integration formula \cite{Xu1999}
for the stiffness matrix, which will be used to prove uniform the boundedness of the scheme:
    \begin{equation}
    \label{edgeint}
    \begin{split}
    \mathcal{A}_h(\gamma; \alpha, \beta):=&\;
 \sum_{K\in\Omega_h} \sum_{E\in \mathcal{E}_K} \omega^K_{E} \Tilde\gamma_{E}\delta_{E}(\alpha)\delta_{E}(\beta)\\
  =&\; \sum_{E\in\mathcal{E}_h} \omega_E
  {\Tilde\gamma}_{E}\delta_{E}(u_h)\delta_{E}(v_h),
    \end{split}
    \end{equation}
    where 
    $ \omega_E:=\sum_{K \supset E} \omega^K_{E}$.
    Here for an edge $E$ with vertices $v_i$
    and $v_j$, we have
    \begin{equation}
        \delta_{E}(u_h) = u_h(v_i) - u_h(v_j).
    \end{equation}
    The quantity $\Tilde\gamma_E$ is the following harmonic average on $E$,
    \begin{equation}
        \Tilde\gamma_E = \left[\frac{1}{|E|}\int_E \frac{1}{\gamma} \,ds \right]^{-1}.
    \end{equation}
    Also, the weights $\omega^K_E$ are given by the identity \cite{Xu1999,barth1992aspects}:
    \begin{equation}
    \omega_E^K=\frac{1}{d(d-1)}\left|\kappa_E^K\right| \cot \theta_E^K,
    \end{equation}
    with $d$ being the number of dimensions, $\theta_E^K$ being the angle between faces not containing edge $E$, and $\kappa_E^K$ the ($d-2$) dimensional simplex formed by their intersection.

\subsection{Properties}
In this section, we prove several important results for the scheme \eqref{semi-imp}. We begin by showing that the scheme satisfies discrete versions of the properties in \eqref{prop}.
\begin{theorem}
\label{thm:pp}
    The fully discrete, semi-implicit scheme \eqref{semi-imp} conserves mass, preserves positivity, and dissipates energy in the following form
    \begin{equation}
        E_h^\n - E_h^\nm \leq -\int_{\Omega_h} m \Delta t \exp \left(m u_h^\nm\right) |\nabla u_h^\n|^2 dx,
    \end{equation}
    where the energy 
    $
    E_h^n:=\mathcal{M}(\exp(u_h^n), u_h^n-1).
    $
\end{theorem}
\begin{proof}
     The positivity of the scheme is guaranteed by the log-density variable formulation since $\exp(u_h^\n) >0$ always. The mass conservation can be proved by picking $v_h=1$ in \eqref{semi-imp}:
     \begin{equation}
         \label{massproof}
         \mathcal{M}\left(\frac{\exp (u_{h}^\n)-\exp (u_h^\nm)}{\Delta t}, 1\right) + 0 = 0 \implies \mathcal{M}\left(\exp (u_{h}^\n), 1\right) = \mathcal{M}\left(\exp (u_{h}^\nm), 1\right).
     \end{equation}
The only non-trivial property to prove is energy dissipation. Observe that by Taylor expansion, we get
\begin{equation}
\label{Taylor1}
(\exp(a)-\exp(b)) a= \exp(a)(a-1)- \exp(b)(b-1)+\frac{1}{2} \exp (\xi)(a-b)^2,
\end{equation}
where $\xi$ is a function between $a$ and $b$. Now, picking $v=u_h^\n$ and using the above Taylor expansion, we get
\begin{equation}
\label{energyproof}
\begin{split}
    \mathcal{M} \left(\exp (u_{h}^\n), u_h^\n-1\right) - \mathcal{M} \left(\exp (u_{h}^\nm), u_h^\nm-1\right) = &-\Delta t 
    \mathcal A \left( m\exp (m u_h^\nm), \nabla u_h^\n, \nabla u_h^\n \right); \\
    & \quad - \mathcal M \left( \exp(\xi), \frac{(u_h^\n-u_h^\nm)^2}{2} \right)
\end{split}
\end{equation}
where $\xi$ is a function between $u_h^\n$ and $u_h^\nm$.
This completes the proof.
\end{proof}

Next we prove unique solvability of the scheme 
\eqref{semi-imp} where mass lumping \eqref{mlquad}
is used for the mass operator.
\begin{theorem}
    The fully discrete scheme \eqref{semi-imp} is uniquely solvable provided that mass lumping \eqref{mlquad} is used to evaluate the mass operator.
\end{theorem}

\begin{proof}
We prove this result using matrix-vector notation. 
Denoting $\underline{u}^j$ as the coefficient vector of solution $u_h^j\in V_h$, 
the scheme \eqref{semi-imp} with mass lumping
can then be written in the following matrix-vector form:
    \begin{equation}
    \label{mat-vec}
        \bld M (\exp(\underline{u}^\n)-\exp(\underline{u}^\nm)) + \Delta t \bld A^\nm \underline{u}^\n = 0,
    \end{equation}
    where $\bld M$ is the diagonal mass matrix, by virtue of mass-lumping, with entries
    $\bld M_{ii} = \frac{|S_i|}{(d+1)},$ and $\bld A^\nm$ is a symmetric positive semidefinite stiffness matrix with entries:
    \begin{equation}
        \bld A^\nm_{ij} =
        \mathcal{A}\left(m\exp (m u_h^\nm); 
\nabla \phi_i, \nabla \phi_j\right).
    \end{equation}

It is clear that the nonlinear system \eqref{mat-vec} is the Euler-Lagrange equation of the following
minimization problem:
\[
\underline u^n:=\mathrm{argmin}_{\underline u\in \mathbb{R}^{N_V}}         F(\underline u), 
\]
where the energy 
functional is
       \begin{equation}
    \label{functional}
         F(\underline u) 
         = \underline{1}\cdot\bld M\exp(\underline u) - \underline u\cdot \bld M\exp(\underline u^\nm)+ \frac{\Delta t}{2}\underline u \cdot \bld A^\nm \underline u,
 \end{equation}
 where $\underline{1}$ is the vector of ones of size $N_V$.
 Hence, unique solvability of \eqref{mat-vec}
 is equivalent to the coercivity 
 (existence) and strictly convexity (uniqueness) of the functional \eqref{functional}.
Both properties can be easily verified using the positivity of the mass and stiffness matrices and elementary calculation. We leave out the details.

\end{proof} 

Finally, we prove the uniform boundedness of the scheme \eqref{semi-imp} when 
mass lumping \eqref{mlquad} is used for the
mass operator, and the edge-based integration
 \eqref{edgeint}
is used for the stiffness operator.
\begin{theorem}
\label{thm:bdd}
If the triangulation $\Omega_h$ is Delaunay, then the solution to the fully discrete scheme \eqref{semi-imp} with mass lumping \eqref{mlquad} being used for the mass operator and edge-based integration \eqref{edgeint} being used for the stiffness operator is uniformly bounded.
That is, given $0<\varepsilon_1<\varepsilon_2$ such that
    $\varepsilon_1 \leq \exp(u_h^\nm) \leq \varepsilon_2$, we have $\varepsilon_1 \leq \exp(u_h^\n) \leq \varepsilon_2$. 
\end{theorem}

\begin{proof}
 We only prove the lower bound, i.e., if 
 $\exp(u_h^\nm)\ge\varepsilon_1$, then $\exp(u_h^n)\ge\varepsilon_1$, as the upper bound use the same argument.

 Taking a non-negative test function $v_h\in V_h$ in \eqref{semi-imp} such that its coefficient $v_i=\max\{\epsilon_1-\exp(u_i^n), 0\}$, and using the matrix-vector form \eqref{mat-vec}, we get 
 \begin{align}
       \label{test1}
       \underline v\cdot \bld M (\exp(\underline{u}^\n)-\exp(\underline{u}^\nm)) + \Delta t  \underline v\cdot\bld A^\nm \underline{u}^\n = 0.
 \end{align}
Using the fact the $\bld M$ is a diagonal positive matrix, $\exp(\underline{u}^\nm)-\varepsilon_1\ge 0$, and definition of $\underline v$,  
we have 
 \[
 \underline v\cdot \bld M (\exp(\underline{u}^\n)-\exp(\underline{u}^\nm)) 
 =
 \underbrace{\underline v\cdot \bld M (\exp(\underline{u}^\n)-\varepsilon_1)}_{\le 0}-
\underbrace{\underline v\cdot \bld M (\exp(\underline{u}^\nm))-\varepsilon_1) }_{\ge 0}
\le 0.
 \]
 Next, using the edge integration formula \eqref{edgeint}, we have 
 \begin{align}
 \label{test2}
 \underline v\cdot\bld A^\nm \underline{u}^\n
 = \sum_{E_{ij}\in\mathcal{E}_h}\omega_{E_{ij}}\Tilde{\gamma}_{E_{ij}} (v_i-v_j)(u_i-u_j),
 \end{align}
 where $\gamma=m\exp(u_h^\nm)$.
 By definition of $v_h$, we have 
\[
(v_i-v_j)(u_i-u_j)
=
(v_i-v_j)((u_i-\varepsilon_1) - (u_j-\varepsilon_1))
\le 0.
\]
Hence, the term \eqref{test2} is non-positive as long as $\omega_{E_{ij}}>0$ for all $E_{ij}\in\mathcal{E}_h$, which is equivalent to the requirement that the triangulation $\Omega_h$ is Delaunay; see \cite{Xu1999}.
In this case, we have 
\[
\underline v\cdot \bld M (\exp(\underline{u}^\n)-\varepsilon_1) = 
\sum_{i=1}^{N_V}
\bld M_{ii}(\exp({u}_i^n)-\varepsilon_1)
\max\{\epsilon_1-\exp(u_i^n), 0\} = 0
\]
thanks to \eqref{test1}. 
Hence, $\exp(u_i^n)\ge \varepsilon_1$, which completes the proof.

\end{proof}

\begin{remark}
We remark that the edge-based integration 
in Theorem \ref{thm:bdd} above is only used to prove the non-positivity of the term \eqref{test2} on Delaunay triangulations. When the mesh is a structured tensor-product grid, such bound preservation can be easily proven with other standard numerical integration rules; see, e.g., \cite{GU2020109378}.

Moreover, although we use the edge integration formula \eqref{edgeint} to prove the uniform boundedness result in Theorem \ref{thm:bdd}, in practical implementation, we simply use the mass-lumping quadrature \eqref{mlquad} to compute both $\bld M$ and $\bld A^\nm$. The resulting scheme is still very robust, and is provable unique solvable and satisfy the properties in Theorem \ref{thm:pp}.
\end{remark}

\begin{remark}
We conclude this section by remarking on the
practical implementation of the scheme \eqref{semi-imp} with compactly supported initial data $\rho^0$. 

In this case, we still interpolate the initial data using \eqref{ic}. Note that we have $u_i^0=-\infty$
whenever $\rho^0(v_i)=0$. Now, the $k$-th Newton iteration for \eqref{mat-vec} takes the form
\begin{equation}
    \left(\bld{M} \bld{D}^{(k-1)} + \Delta t \bld A^\nm \right)\underline{u}^{(k)} = \bld M \left( \bld D^{(k-1)} \underline u^{(k-1)} -\exp{(\underline u^{(k-1)})} + \exp{(\underline u^{n-1})}   \right),
\end{equation}
where $\bld D^{(k-1)}$ is a diagonal matrix with diagonal entry $\bld D^{(k-1)}_{ii} = \exp{u_i^{(k-1)}}$, and $\underline u^{(0)} = \underline u^\nm$. To solve the above linear system for
$\underline u^{(k)}$, we only activate the $i$-th degree of freedom if the diagonal entry $\left(\bld{M} \bld{D}^{(k-1)} + \Delta t \bld A^\nm \right)_{ii}$ is greater than a small cutoff value (e.g., $10^{-14}$), and set the inactive degrees of freedom to $-\infty$.
\end{remark}

\section{Mixed Method}
In this section, we develop our mixed formulation to solve the PME \eqref{PME}.

The physical energy for the PME is given by $U(\rho) = \frac{\rho^m}{m-1}$. We define a new potential variable equal to the derivative of the physical energy:
\begin{equation}
    \mu = U'(\rho) = \frac{m}{m-1}\rho^{m-1}.
\end{equation}
We also define a velocity variable, $\bld u$, and set it equal to the negative gradient of the potential, i.e.,
\begin{equation}
    \bld u = -\nabla \mu.
\end{equation}
Notice that using this definition, we easily observe that
\begin{equation}
    \nabla \rho ^m = \rho\nabla \mu = -\rho\bld u.
\end{equation}
We have now successfully reformulated the PME \eqref{PME} into the
following first-order system
\begin{equation}
\label{pme-mix}
    \begin{cases} 
    &\rho_t+\nabla \cdot (\rho \bld u) =0, \\
    & \bld u+\nabla \mu = 0, \\ 
    & \mu = \frac{m}{m-1}\rho^{m-1},
    \end{cases}
\end{equation}
where the unknowns are density, potential, and velocity. 
We again equip the PME system \eqref{pme-mix} with the homogeneous Neumann boundary condition 
\[
\bld u\cdot\bld n=0 \text{ on }\partial \Omega.
\]
\subsection{Spatial and Temporal Discretizations}
Our energy-stable mixed method can be constructed on structured tensor product meshes in any space dimension, or on general unstructured Delaunay triangular meshes in two dimensions. 
Due to the use of velocity mass lumping as a key tool to establish
the energy stability result, the method does not work on general 
simplicial meshes in three dimensions.
Below we describe the method in detail on 2D triangular meshes, following the meshing notation as in Section \ref{sec:notation}.
For any triangular element $K\in\Omega_h$, we denote 
$\partial K$ as its boundary whose associated outward unit normal is $\bld n_K$.

We shall use the following two finite element spaces:
\begin{equation}
\begin{aligned}
&Q_{h}=\left\{q_{h} \in L^{2}(\Omega) : \enspace  \left.q_{h}\right|_{K} \in P_{0}(K), \enspace \forall K \in \Omega_{h}\right\}, \\
&\bld V_{h}=\left\{\bld v_{h} \in \operatorname{H}(\operatorname{div}; \Omega) : \enspace \left.\bld v_{h}\right|_{K} \in RT_0(K), \enspace \forall K \in \Omega_{h}, \;\;
\bld v_h\cdot\bld n = 0,\quad\text{ on } \partial\Omega
\right\},
\end{aligned}
\end{equation}
where $P_0(K)$ is the constant space, and $RT_0(K)=[P_0(K)]^2\oplus \bld xP_0(K)$ is the lowest-order Raviart-Thomas space on the triangle $K$. Note that the Neumann boundary condition is encoded in the velocity space $\bld V_h$.

The spatial discretization now reads as follows: find $\rho_h, \mu_h \in Q_h$ and $\bld u_h \in \bld V_h$ such that, for $t>0$,
\begin{subequations}
\label{m:space}
\begin{align}
\label{m:sd2}
\sum_{K \in \Omega_h}\left[ \int_K (\rho_h)_t q_h\,\mathrm{dx} + \int_{\partial K} \hat{\rho}_h\bld u_h \cdot \bld n_K q_h\,\mathrm{ds} \right]
=&\; 0 , \qquad \forall \enspace q_h \in Q_h,\\
    \label{m:sd3}
    \sum_{K \in \Omega_h}\left[ \int_K \bld u_h \cdot \bld v_h\,\mathrm{dx} - \int_{\partial K}  \mu_h \bld v_h \cdot \bld n_K\,\mathrm{ds} \right ] = &\;0, \qquad \forall \enspace \bld v_h \in \bld V_h,\\
\label{m:sd1}
    \sum_{K \in \Omega_h} \int_K (\mu_h-\tfrac{m}{m-1} \rho_h^{m-1}) r_h dx = &\;0,  \qquad \forall \enspace r_h \in Q_h,
\end{align}
\end{subequations}
where
$\hat \rho_h$ 
is the upwinding numerical flux, i.e., 
 
given an edge $E$ shared by two elements $K^+$ and $K^-$ with $\bld u_h\cdot\bld n_{K^-}|_E\ge0$, $\hat{\rho}_h$ takes 
value from $K^-$:
\begin{equation}
\left. \hat{\rho}_h \right|_{E}
= (\rho_h|_{K^-})|_{E}, \quad \text{ where }
\bld u_h\cdot\bld n_{K^-}|_E\ge0.
\end{equation} 

Analogous to the log-density approach, we apply a first-order semi-implicit discretization of the ODE system \eqref{m:space} to obtain the following fully discrete scheme: Given time step size $\Delta t$ and data $\rho_h^\nm \in Q_h$ at time $t^\nm$, find $\mu^n_h, \rho^n_h \in Q_h$, and $\bld u_h^n \in \bld V_h$, at time $t^\n$, such that

\begin{subequations} \label{m:fd}
\begin{align}
\mathcal M_h \left( \frac{\rho_h^n-\rho_h^{n-1}}{\Delta t}, q_h \right) + \mathcal A_h\left(\hat \rho_h^\nm ; \bld u_h^n, q_h  \right) = 0, \quad &\forall q_h \in Q_h, \label{m:fd2}\\
\bar{\mathcal M_h}(\bld u_h^n, \bld v_h) - \bar{\mathcal{A}_h} (\mu^n_h, \bld v_h) = 0, \quad &\forall \bld v_h \in \bld V_h, \label{m:fd3}\\
\mathcal{M}_h \left(\mu^n_h-\frac{m}{m-1}(\rho_h^\n)^{m-1}
,r_h \right) 
= 0, \quad &\forall r_h \in Q_h, \label{m:fd1}
\end{align}
\end{subequations}
where the associated operators are defined as:
\begin{subequations}
    \begin{align}
        &\mathcal{M}_h \left(\alpha,\beta \right) = \sum_{K \in \Omega_h} \int_K \alpha \beta  \mathrm{dx}, 
        && \bar{\mathcal M}_h(\boldsymbol \alpha, \boldsymbol \beta) = \sum_{K \in \Omega_h} \int_K \boldsymbol\alpha \cdot \boldsymbol\beta \mathrm{dx},\\
        &\mathcal A_h\left(\gamma ; \boldsymbol \alpha, \beta  \right) = \sum_{K \in \Omega_h} \int_{\partial K} \gamma \boldsymbol\alpha \cdot \bld n_K \beta  \mathrm{dx}, 
        && \bar{\mathcal A}_h\left( \alpha, \boldsymbol\beta  \right) =  \sum_{K \in \Omega_h} \int_{\partial K} \alpha \boldsymbol \beta \cdot \bld n_K  \mathrm{dx}.
    \end{align}
\end{subequations}
Again, the numerical flux $\hat \rho_h^\nm$ in \eqref{m:fd2} is taken to be the upwind flux:
given an edge $E$ shared by two elements $K^+$ and $K^-$ with $\bld u_h^n\cdot\bld n_{K^-}|_E\ge0$, 
we take
\begin{equation}
\label{m:flux}
\left. \hat{\rho}_h^{n-1} \right|_{E}
= (\rho_h^{n-1}|_{K^-})|_{E}, \quad \text{ where }
\bld u_h^n\cdot\bld n_{K^-}|_E\ge0,
\end{equation}
Standard one-point numerical integration rules are used in the evaluation of the operators $\mathcal{M}_h, \mathcal{A}_h,$ and $\bar{\mathcal A}_h$. However, the energy stability of the scheme requires that the velocity mass matrix be mass lumped. As such, the operator $\bar{\mathcal M}_h$ must be evaluated using an appropriate mass-lumping quadrature. 
On tensor-product grids, this is readily achieved by the use of the trapezoidal rule. On triangular meshes, mass-lumping is achieved by 
using the following formula given in \cite{Baranger1996}:

\begin{equation}
    \label{m:mlquad}
    \bar{\mathcal M_h }\left( \boldsymbol \alpha,\boldsymbol \beta \right) := \sum_{K \in \Omega_h} \sum_{E \in \mathcal E_K}  \omega_E^K \varphi_E(\boldsymbol \alpha)\varphi_E(\boldsymbol \beta) = \sum_{E \in \mathcal{E}_h} \omega_E\varphi_E(\boldsymbol \alpha)\varphi_E(\boldsymbol \beta),
\end{equation}
where $\omega_E := \sum_{K \supset E} \omega_E^K$, and $\varphi_E(\boldsymbol \alpha):=\boldsymbol \alpha\cdot\bld n_E$ denotes the normal flux of $\boldsymbol\alpha$ through edge $E$. The weights $\omega_E^K$ are given by
\begin{equation}
    \omega_E^K=\frac{1}{2} \cot \theta_E^K,
\end{equation}
where $\theta_E^K$ is the angle opposite to edge $E$ in $K$. 
It is clear that the mass matrix associated with the integration rule in \eqref{m:mlquad} is a diagonal matrix, whose diagonal entries are positive provided that the mesh is Delaunay.

To efficiently solve the scheme \eqref{m:fd}, 
we first apply static condensation to locally solve 
the potential and velocity variables 
$\mu_h^n$ and $\bld u_h^n$
as functions of the density $\rho_h^n$ using 
equations \eqref{m:fd3}--\eqref{m:fd1}, and then 
use Newton's method to solve the resulting parabolic nonlinear system \eqref{m:fd2} for density alone.

\subsection{Properties}
\begin{theorem}
    \label{m:massenergythm}
    Provided a Delaunay triangulation $\Omega_h$, the fully discrete scheme \eqref{m:fd} is mass conservative and energy dissipative in the following forms
    \begin{subequations}
        \begin{align}
            \mathcal M_h \left( \rho_h^n, 1 \right) &= \mathcal M_h \left( \rho_h^\nm, 1 \right) 
            \label{m:masscon}\\
             \mathcal M_h \left( U(\rho^n),1 \right) - \mathcal M_h \left( U(\rho^\nm),1 \right) &\leq -\Delta t \sum_{E \in \mathcal{E}_h} \omega_E \varphi_E(\bld \alpha) \varphi_E(\bld u_h^n)^2 \label{m:energydiss},
        \end{align}
    \end{subequations}
    where $U^n = \frac{m}{m-1}(\rho^n)^{m-1}$ is the physical energy, and the right-hand side in \eqref{m:energydiss} comes from the quadrature formula \eqref{m:mlquad}.
\end{theorem}
\begin{proof}
    By picking $q_h = 1$ in \eqref{m:fd2} and using the homegenous boundary condition $\bld u_h^n\cdot\bld n|_{\partial\Omega}=0$, we can easily recover the mass conservation \eqref{m:masscon}. Now, to prove energy stability, we take the test function 
    $\bld v_h$ in \eqref{m:fd3} to be a function in $\bld V_h$ such that its normal flux through edge $E$ is given by $\varphi_E(\bld v_h) = \varphi_E(\hat\rho_h^\nm\bld u_h^n)$. For this choice of $\bld v_h$, we have
    \begin{equation}
    \label{m:Aequality}
        \mathcal A_h\left(\hat \rho_h^\nm ; \bld u_h^n, \mu^n_h  \right) = \bar{\mathcal{A}_h} (\mu^n_h, \bld v_h).
    \end{equation}
    Additionally, through an application of \eqref{m:mlquad}, we obtain
    \begin{equation}
        \label{m:velmassml}
        \bar{\mathcal M_h}(\bld u_h, \bld v_h) = \sum_{E \in \mathcal E_h} \omega_E\hat\rho_h^\nm \varphi_E(\bld u_h^n)^2
    \end{equation}
Next taking test function 
$r_h=\frac{\rho_h^n-\rho_h^\nm}{\Delta t}$
in \eqref{m:fd1}, we get 
\begin{equation}
\label{m:taylorexp}
    \begin{split}
        &
        \mathcal M_h \left(\frac{\rho_h^n-\rho_h^{n-1}}{\Delta t}, \mu_h^n \right) 
        = \mathcal M_h \left(\frac{\rho_h^n-\rho_h^{n-1}}{\Delta t}, U'(\rho_h^n) \right)\\
        &= \frac{1}{\Delta t}\mathcal M_h \left( U(\rho_h^n),1 \right) - \frac{1}{\Delta t}\mathcal M_h \left( U(\rho_h^\nm),1 \right) 
        + \frac{1}{\Delta t} \mathcal M_h\left( \frac{U''(\xi)}{2}(\rho_h^n-\rho_h^\nm)^2 \right) ,
    \end{split}
\end{equation}
where we have used Taylor expansion with $\xi$ being a function between $\rho^n_h$ and $\rho_h^\nm$. 
Finally, taking the test function 
$q_h=\mu_h^n$ in \eqref{m:fd2}, and 
using the above relations, we get
\begin{equation}
    \begin{split}
        \mathcal M_h \left( U(\rho^n),1 \right) - \mathcal M_h \left( U(\rho^\nm),1 \right) &= -\Delta t\sum_{E \in \mathcal E_h} \omega_E\hat\rho_h^\nm \varphi_E(\bld u_h^n)^2 -  \mathcal M_h\left( \frac{U''(\xi)}{2}(\rho_h^n-\rho_h^\nm)^2 \right)\\
    &\leq -\Delta t\sum_{E \in \mathcal E_h} \omega_E\hat\rho_h^\nm \varphi_E(\bld u_h^n)^2
    \end{split}
\end{equation}
The right-hand side is non-positive if $\omega_E \geq 0$ which is guaranteed if the triangulation $\Omega_h$ is Delaunay. This completes the proof.
\end{proof}

Next, we prove that the scheme \eqref{m:fd} is positivity preserving under a usual CFL time stepping constraint.

\begin{theorem}
   Given $\rho^\nm_h \geq 0$, the fully discrete scheme \eqref{m:fd} is positivity preserving under the CFL condition,
   \begin{equation}
   \label{cfl}
{\Delta t} \sum_{E \in \mathcal E^-_K} |\bld u_h^n\cdot{\bld n_K}|_E\frac{|E|}{|K|} \leq 1, \quad \forall K \in \Omega_h,
   \end{equation}
   where 
\[
\mathcal{E}_K^-:=\{E\in\mathcal{E}_K: \quad
\bld u_h^n\cdot\bld n_K|_E\ge 0.
\}
\]
\end{theorem}
\begin{proof}
Restricting the mass conservation equation \eqref{m:fd2} to a single element $K\in \Omega_h$, we have 
    \begin{equation}
        \label{m:Keqn}
        \frac{|K|}{\Delta t}(\rho_K^n - \rho_K^\nm) + \sum_{E \in \mathcal E_K} \hat\rho^\nm_E \left.(\bld u^n_h \cdot \bld n_K)\right|_E |E| = 0,
    \end{equation}
    where $\rho_K^n:=\rho_h^n|_K$ is the restriction of the function $\rho_h^n\in Q_h$ to the element $K$.
        This implies that 
    \begin{align*}
\rho_K^n = \rho_K^\nm -\frac{\Delta t}{|K|} \sum_{E \in \mathcal E_K} \hat\rho^\nm_E\left.(\bld u^n_h \cdot \bld n_K)\right|_E |E|,
    \end{align*}   
By the definition of
$\mathcal{E}^-_K$, we have 
$\hat\rho_E^{n-1}=\rho_K^{n-1}$ for all $E\in\mathcal{E}^-_K$.
Hence, 
    \begin{align*}
\rho_K^n = \rho_K^\nm (1-\frac{\Delta t}{|K|}
\sum_{E \in \mathcal E^-_K} \left|\bld u^n_h \cdot \bld n_K\right|_E  |E|)
+
\frac{\Delta t}{|K|}
\sum_{E \in \mathcal E_K\backslash\mathcal{E}^-_K} \hat\rho_E^{n-1}\left|\bld u^n_h \cdot \bld n_K\right|_E  |E|,
    \end{align*}
    Both terms on the right hand side are nonnegative 
    under the assumption that $\rho_h^{n-1}\ge0$ and 
    \eqref{cfl}.
    This completes the proof.
\end{proof}

\begin{remark}
We note that for tensor product meshes, the stronger property of unconditional positivity preservation 
can be proven following similar arguments as in the proof of Theorem \ref{thm:bdd}; see also \cite{GU2020109378}. 
The key is to locally eliminate velocity and potential degrees of freedom to express the scheme as a  finite volume scheme for the piecewise constant density unknown only. 
In this case, unconditionally positivity preservation holds for any consistent numerical flux, that is, the upwinding flux is not needed for the positivity proof on tensor-product meshes. We leave out the detailed derivation.
\end{remark}

\section{Numerical Results}
\label{sec:num}
In this section, we present out numerical findings. All computations are performed using the Python interface of the open-source library NGSolve \cite{Schoberl16}. The source code for these computations is available at the following git repository: \url{https://github.com/avj-jpg/pme}.
\begin{table}[]
\resizebox{\textwidth}{!}{
\begin{tabular}{|ccccccc|ccccc|}
\hline
  &                          & \multicolumn{5}{|c|}{Log-density method}                                                        & \multicolumn{5}{c|}{Mixed method}                                                                                                                       \\ \cline{3-12} 
m & \multicolumn{1}{c|}{N}   & $\Delta t$    & \begin{tabular}[c]{@{}c@{}}Error in\\ $[-5,5]$\end{tabular} & Order & \begin{tabular}[c]{@{}c@{}}Error in \\ $[-10,10]$\end{tabular}     & Order & $\Delta t$   & \begin{tabular}[c]{@{}c@{}}Error in \\ $[-5,5]$\end{tabular} & Order & \begin{tabular}[c]{@{}c@{}}Error in \\ $[-10,10]$\end{tabular} & Order \\ \hline
2 & \multicolumn{1}{c|}{100} & 1/5   & 1.19e-01                                                    & -     & 3.72e-01 &  -     & 1/10 & 4.53e-02                                                     &  -     & 8.48e-02                                                            &  -     \\
  & \multicolumn{1}{c|}{200} & 1/20  & 3.04e-02                                                    & 1.971 & 1.03e-01 & 1.860 & 1/20 & 2.27e-02                                                     & 0.999 & 4.26e-02                                                         & 0.992 \\
  & \multicolumn{1}{c|}{400} & 1/80  & 7.57e-03                                                    & 2.009 & 2.49e-02 & 2.043 & 1/40 & 1.13e-02                                                     & 0.999 & 2.14e-02                                                         & 0.994 \\
  & \multicolumn{1}{c|}{800} & 1/320 & 1.88e-03                                                    & 2.007 & 6.18e-03 & 2.010 & 1/80 & 5.67e-03                                                     & 1.000 & 1.08e-02                                                         & 0.992 \\ \hline
3 & \multicolumn{1}{c|}{100} & 1/5   & 6.63e-02                                                    & -     & 3.15e-01 & -     & 1/10 & 1.94e-02                                                     & -     & 8.61e-02                                                         & -     \\
  & \multicolumn{1}{c|}{200} & 1/20  & 1.64e-02                                                    & 2.012 & 9.36e-02 & 1.752 & 1/20 & 9.76e-03                                                     & 0.995 & 4.69e-02                                                         & 0.876 \\
  & \multicolumn{1}{c|}{400} & 1/80  & 3.97e-03                                                    & 2.050 & 2.55e-02 & 1.878 & 1/40 & 4.89e-03                                                     & 0.997 & 2.65e-02                                                         & 0.825 \\
  & \multicolumn{1}{c|}{800} & 1/320 & 9.31e-04                                                    & 2.091 & 8.93e-03 & 1.511 & 1/80 & 2.45e-03                                                     & 0.998 & 1.54e-02                                                         & 0.782 \\ \hline
4 & \multicolumn{1}{c|}{100} & 1/5   & 4.07e-02                                                    & -     & 2.33e-01 & -     & 1/10 & 1.19e-02                                                     & -     & 1.05e-01                                                         & -     \\
  & \multicolumn{1}{c|}{200} & 1/20  & 1.02e-02                                                    & 1.994 & 7.51e-02 & 1.634 & 1/20 & 5.95e-03                                                     & 1.006 & 6.03e-02                                                         & 0.793 \\
  & \multicolumn{1}{c|}{400} & 1/80  & 2.47e-03                                                    & 2.050 & 2.70e-02 & 1.478 & 1/40 & 3.01e-03                                                     & 0.980 & 3.25e-02                                                         & 0.892 \\
  & \multicolumn{1}{c|}{800} & 1/320 & 4.43e-04                                                    & 2.479 & 1.24e-02 & 1.126 & 1/80 & 1.49e-03                                                     & 1.016 & 2.14e-02                                                         & 0.606\\
  \hline
\end{tabular}
}\caption{Convergence results for 1D Barenblatt initial data.}
    \label{t:1DBB}
\end{table}
\begin{figure}[h!]
\centering
\begin{subfigure}{0.45\textwidth}
    \includegraphics[width=\textwidth]{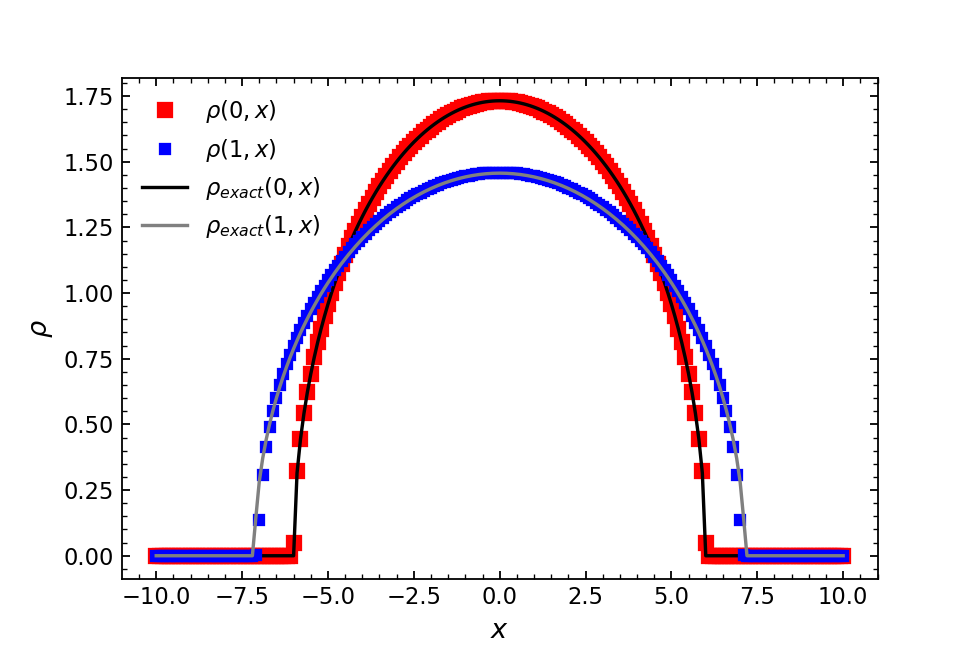}
    \caption{Log density method.}
\end{subfigure}
\begin{subfigure}{0.45\textwidth}
    \includegraphics[width=\textwidth]{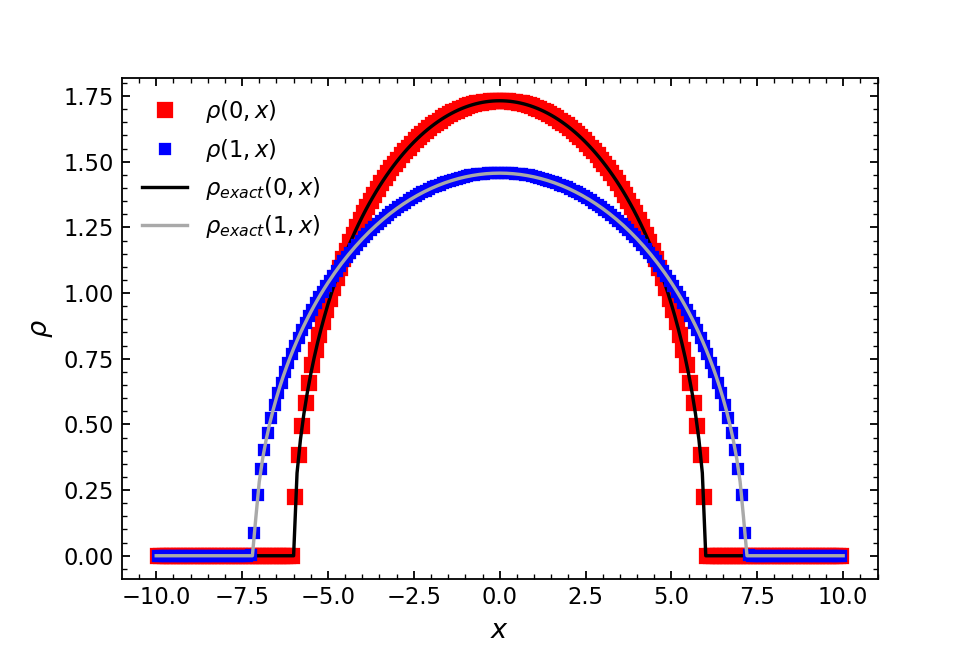}
    \caption{Mixed method.}
\end{subfigure}

\caption{Evolution of 1D Barenblatt initial data by the two schemes for $m=3$, $\Delta t = 0.05$, and $N=200$ elements.}
\label{fig:1DBB}
\end{figure} 
\subsection{1D Barenblatt solution}
The Porous Medium Equation admits an exact weak solution formulated by Barenblatt \cite{Barenblatt1952} and Pattle \cite{10.1093/qjmam/12.4.407}. In the one-dimensional case, the Barenblatt solution is given by the equation: 
\begin{equation}
    \label{1dbb}
    \rho_B(x,t) = (t+1)^{-k} \left(s_0 - \frac{k(m-1)}{2m} \frac{x^2}{(t+1)^{2k}}\right)_+^{\frac{1}{m-1}}, \quad t > 0,
\end{equation}
where $k = (m+1)^{-1}$, and $s_0$ denotes a scaling factor. Note that this data is compactly supported in the interval $[-\eta_m(t), \eta_m(t)]$, where the right boundary $\eta_m(t)$ moves as:
\begin{equation}
    \eta_m(t) = \sqrt{\frac{2m}{k(m-1)}}(t+1)^k.
\end{equation}
To verify the accuracy of the two schemes, we use \eqref{1dbb} as the initial data, with $x 
\in [-10,10]$, $s_0=3$, and $m\in\{2,3,4\}$, and conduct a spacetime mesh-refinement study for the $L^2$-error of density at final time $t=1$. We record the $L^2$-error in the entire domain $[-10,10]$ and in the interval $[-5,5]$ away from the interface. We utilize a sequence of meshes consisting of $\{100 \cdot 2^{i}\}_{i=0}^3$ spatial elements and set the timestep size to $\frac{1}{5\cdot 4^{i}}$  for the log-density method and to $\frac{1}{10\cdot 2^{i}}$ for the mixed method, correspondingly. The results of this convergence study are recorded in table \ref{t:1DBB}. For all three values of $m$, we observe that the log-density method is second-order accurate in space and first-order accurate in time in the region $[-5,5]$ where the solution remains smooth at final time. For larger values of $m$, the order of convergence deteriorates when the $L^2$-error calculated in the entire domain. This is anticipated since the error is expected to be greater at the interface for larger values of $m$, owing to the decreased regularity of the solution. 

On the other hand, we observe that the mixed scheme has the expected first order accuracy in both space and time in the interval $[-5,5]$. Similar to the case of the log-density method, the order of convergence is observed to decay with $m$ if the $L^2$-error is measured in the entire domain $[-10,10]$. However, the decay in order is slower than in the case of the log-density scheme.

Figure \ref{fig:1DBB} displays plots of the numerically computed Barenblatt solution using the two schemes at time $t=1$, for $m=3$ and $\Delta t=0.05$. These plots illustrate that the initial profile is accurately evolved by both schemes, with no oscillations emerging near the interface.

\begin{figure}[h!]
\centering
\begin{subfigure}{0.33\textwidth}
    \includegraphics[width=\textwidth,trim={0 0 3cm 0},clip]{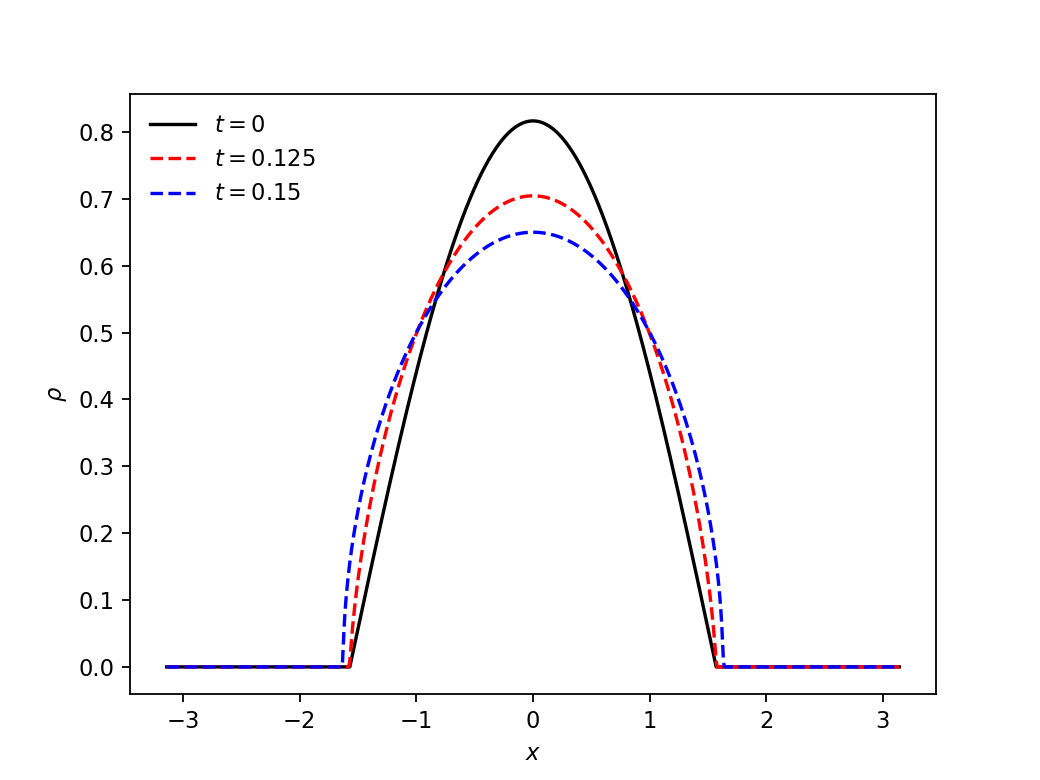}
    \caption{Log-density method}
\end{subfigure}
\begin{subfigure}{0.33\textwidth}
\includegraphics[width=\textwidth,trim={0 0 3cm 0},clip]{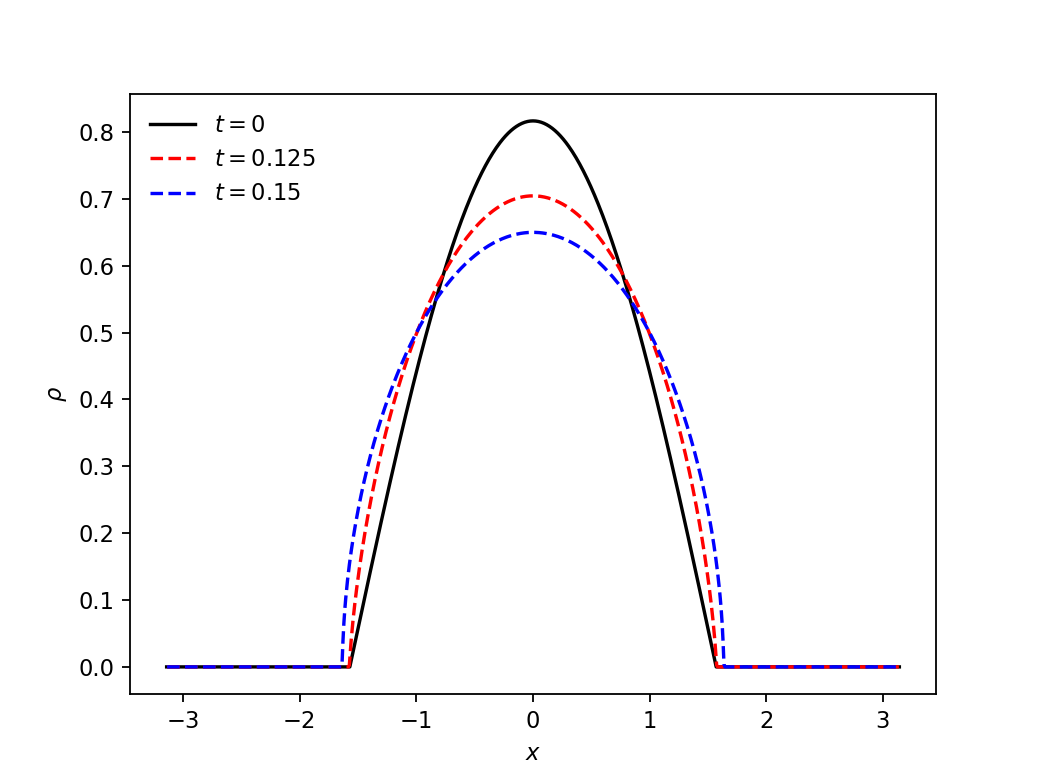}
    \caption{Mixed method}
\end{subfigure}

\begin{subfigure}{0.33\textwidth}
\includegraphics[width=\textwidth,trim={0 0 3cm 0},clip]{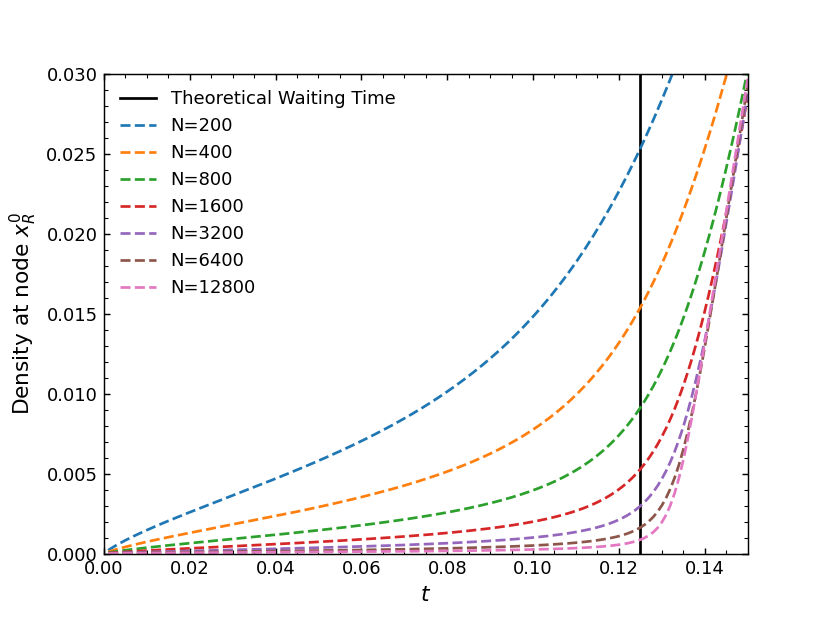}
    \caption{Log-density method}
\end{subfigure}
\begin{subfigure}{0.33\textwidth}
\includegraphics[width=\textwidth,trim={0 0 3cm 0},clip]{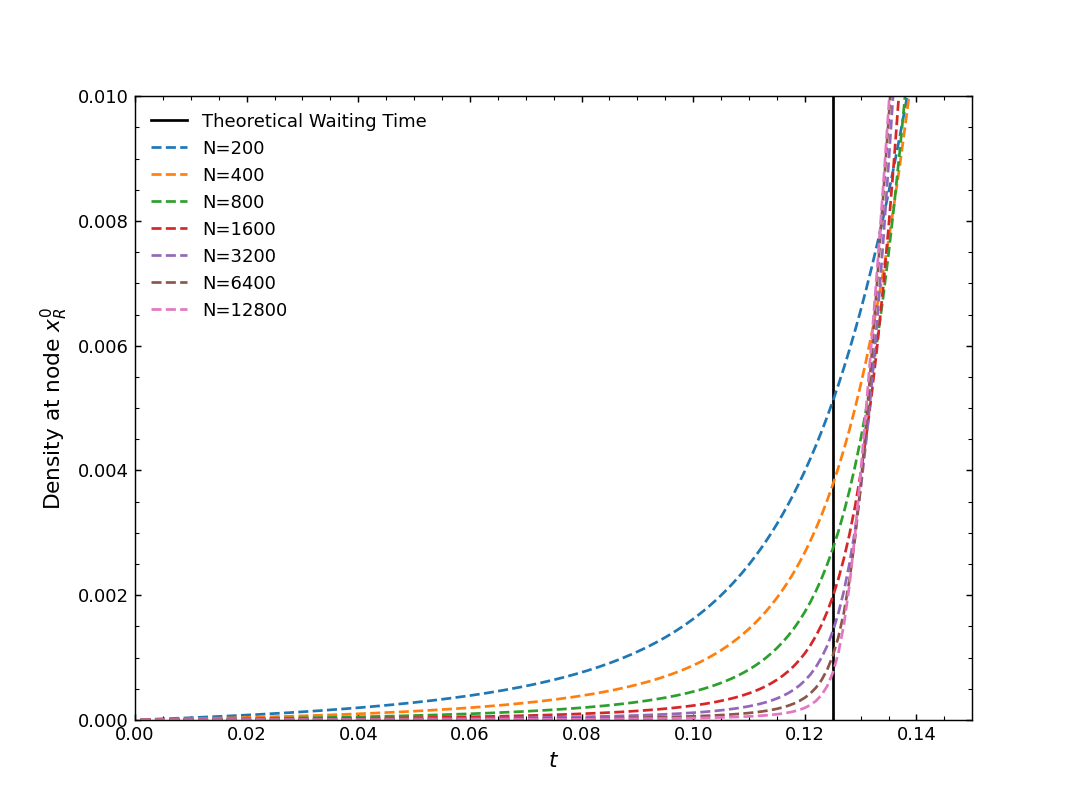}
    \caption{Mixed method}
\end{subfigure}
\caption{Observation of waiting time phenomenon for $m=3$ and $\theta =0$. Panels (A) and (C) show the solutions at various times, and panels (B) and (D) show the density over time at the node ($x_R^0$) corresponding to the initial location of the right interface.}
\label{fig:wt}
\end{figure} 

\subsection{Waiting time phenomenon}
The waiting time phenomenon is a well-known feature of the PME for initial data of form
\begin{equation}
\label{wtdata}
\rho^0(x)= \begin{cases}\left(\frac{m-1}{m}\left((1-\theta) \cos ^2 x+\theta \cos ^4 x\right)\right)^{\frac{1}{m-1}}, & \text { if }-\frac{\pi}{2} \leq x \leq \frac{\pi}{2} \\
0, \quad &\text { otherwise, }\end{cases}
\end{equation}
with $\theta$ in the interval $[0,1]$. For this data, the interface does not move until $t$ is greater than a certain waiting time $t^*$, even though the internal profile continues to evolve. The theoretical waiting time for $0 \leq \theta \leq \frac{1}{4}$ is given by \cite{doi:10.1137/0514049}
\begin{equation}
t^*=\frac{1}{2(m+1)(1-\theta)}.
\end{equation}
In our numerical experiment, we set the parameters $m=3$ and $\theta=0$, and solve the equation \eqref{wtdata} using the two schemes on a sequence of meshes with $\{200\cdot 2^i\}_{i=0}^6$ elements and a fixed time step $\Delta t = 0.001$. The simulation is performed until the final time $t=0.15$, past the theoretical waiting time $t^*=0.125$. We denote by $x_R^0$ the node that corresponds to the location of the right interface at time $t=0$ and track the density value at this node for each time step. 

The results of our experiment are presented in Figure \ref{fig:wt}. Panels (A) and (B) show the density profiles at times $t=0$, $t=0.125$, and $t=0.15$, while (C) and (D) plot the density at $x_R^0$ as a function of time for different values of $N$, as obtained by the log-density and the mixed method, respectively. In both (A) and (B), we observe that the interface remains stationary at $t=0.125$ while the internal profile continues to evolve. At time $t=0.15$, a clear movement of the interface is noticed. The panels (C) and (D) reveal for both methods, when $t$ is equal to the theoretical waiting time (indicated by a black vertical line), the density value at $x_R^0$ gradually decreases towards 0 as the mesh is refined by increasing $N$. Furthermore, for each value of $N$, the mixed method gives a smaller density, and thus a smaller error, at the node $x^R_0$ when $t=t^*$. The smaller error can be ascribed to the piecewise constant spatial discretization of the mixed method, which allows for a sharper capture of the interface, reducing the error near it.

\subsection{Higher dimensional Barenblatt solution}
\begin{table}[]
    \centering
   \resizebox{\textwidth}{!}{
\begin{tabular}{|cc|ccccc|ccccc|}
\hline
  &         & \multicolumn{5}{c|}{Log-density method}                                                                                                            & \multicolumn{5}{c|}{Mixed method}                                                                                                                 \\ \cline{3-12} 
m & N       & dt    & \begin{tabular}[c]{@{}c@{}}Error in \\ $[-3,3]$\end{tabular} & Order & \begin{tabular}[c]{@{}c@{}}Error in\\ $[-6,6]$\end{tabular} & Order & dt   & \begin{tabular}[c]{@{}c@{}}Error \\ in $[-3,3]$\end{tabular} & Order & \begin{tabular}[c]{@{}c@{}}Error in\\ $[-6,6]$\end{tabular} & Order \\ \hline
2 & 32x32   & 1/5   & 2.35e-02                                                     & -     & 6.37e-02                                                    & -     & 1/10 & 1.29e-01                                                     & -     & 2.45e-01                                                    & -     \\
  & 64x64   & 1/20  & 6.61e-03                                                     & 1.828 & 2.97e-02                                                    & 1.101 & 1/20 & 6.46e-02                                                     & 1.004 & 1.25e-01                                                    & 0.971 \\
  & 128x128 & 1/80  & 1.70e-03                                                     & 1.956 & 1.33e-02                                                    & 1.161 & 1/40 & 3.21e-02                                                     & 1.007 & 6.38e-02                                                    & 0.971 \\
  & 256x256 & 1/320 & 4.29e-04                                                     & 1.989 & 5.77e-03                                                    & 1.202 & 1/80 & 1.61e-02                                                     & 0.999 & 3.25e-02                                                    & 0.973 \\ \hline
3 & 32x32   & 1/5   & 9.55e-03                                                     & -     & 1.98e-01                                                    & -     & 1/10 & 7.47e-02                                                     & -     & 3.77e-01                                                    & -     \\
  & 64x64   & 1/20  & 2.99e-03                                                     & 1.677 & 1.21e-01                                                    & 0.706 & 1/20 & 3.72e-02                                                     & 1.006 & 2.01e-01                                                    & 0.905 \\
  & 128x128 & 1/80  & 7.81e-04                                                     & 1.936 & 6.98e-02                                                    & 0.795 & 1/40 & 1.85e-02                                                     & 1.011 & 1.13e-01                                                    & 0.835 \\
  & 256x256 & 1/320 & 1.99e-04                                                     & 1.970 & 4.01e-02                                                    & 0.801 & 1/80 & 9.24e-03                                                     & 0.999 & 6.22e-02                                                    & 0.859 \\ \hline
4 & 32x32   & 1/5   & 5.97e-03                                                     & -     & 3.40e-01                                                    & -     & 1/10 & 4.39e-02                                                     & -     & 5.28e-01                                                    & -     \\
  & 64x64   & 1/20  & 1.55e-03                                                     & 1.949 & 2.33e-01                                                    & 0.547 & 1/20 & 2.19e-02                                                     & 1.002 & 2.86e-01                                                    & 0.884 \\
  & 128x128 & 1/80  & 4.04e-04                                                     & 1.936 & 1.46e-01                                                    & 0.672 & 1/40 & 1.09e-02                                                     & 1.009 & 1.72e-01                                                    & 0.729 \\
  & 256x256 & 1/320 & 1.03e-04                                                     & 1.977 & 8.92e-02                                                    & 0.712 & 1/80 & 5.45e-03                                                     & 0.998 & 1.02e-01                                                    & 0.764 \\ \hline
\end{tabular}
}
    \caption{Convergence results for 2D Barenblatt initial data.}
    \label{t:2DBB}
\end{table}
\begin{figure}[h!]
\begin{subfigure}{.24\textwidth}
\centering
\includegraphics[width=\linewidth,trim={0 0 5.5cm 0},clip]{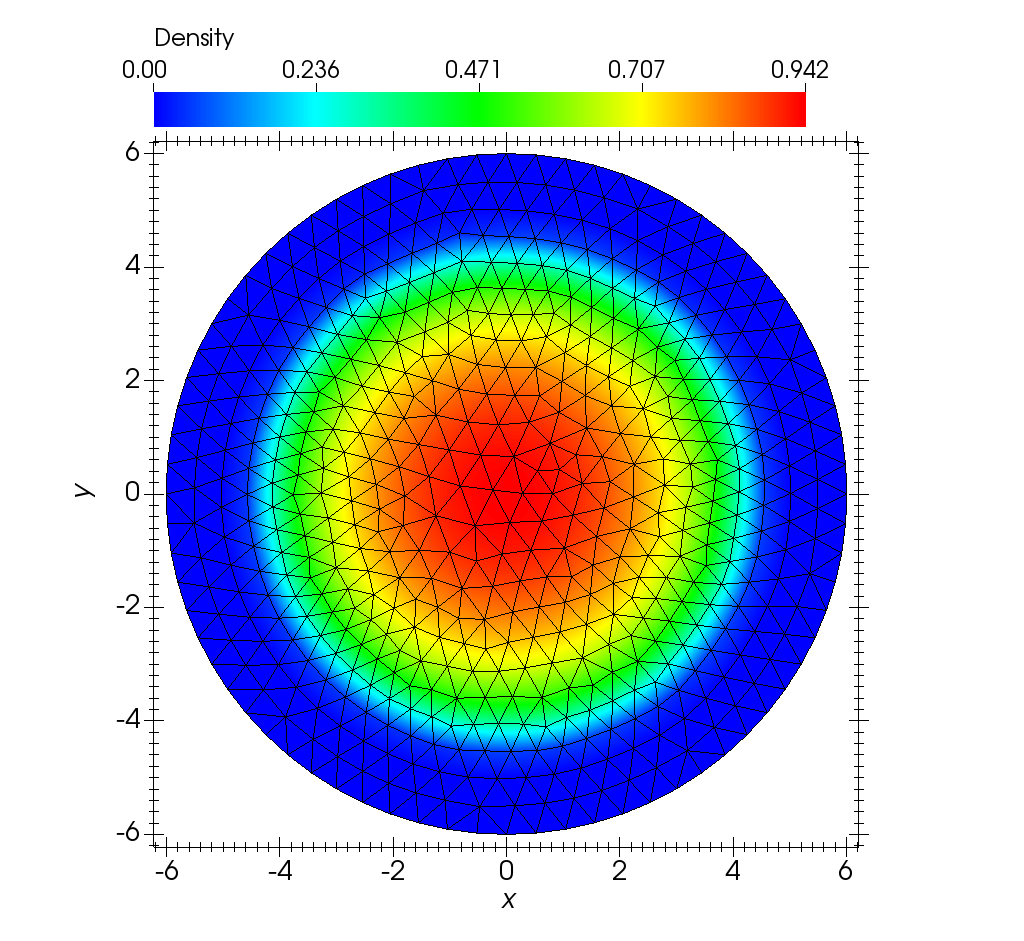}
\centering
\captionsetup{justification=centering}\caption{Density; triangular elements}
\end{subfigure}
\begin{subfigure}{.24\textwidth}
\centering
\includegraphics[width=\linewidth,trim={0 0 5.5cm 0},clip]{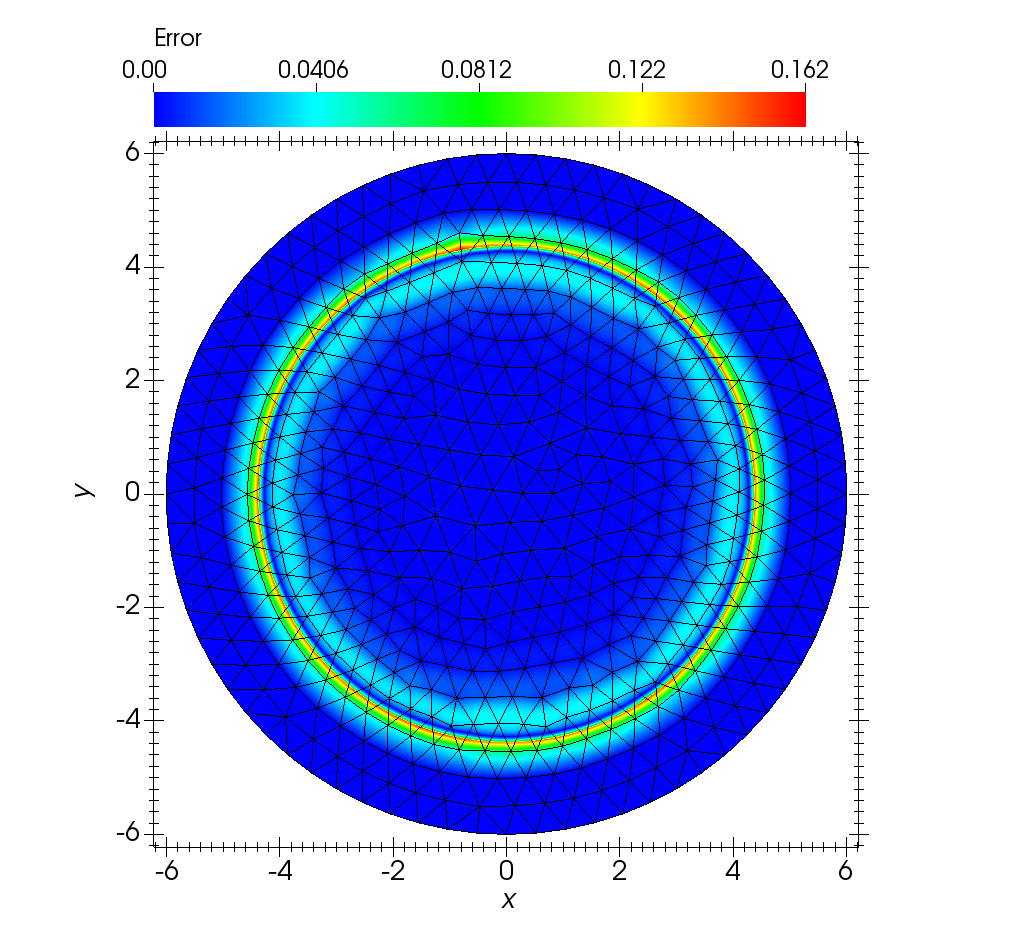}
\captionsetup{justification=centering}\caption{Error; triangular elements}
\end{subfigure}
\begin{subfigure}{.24\textwidth}
\centering
\includegraphics[width=\linewidth,trim={0 0 5.5cm 0},clip]{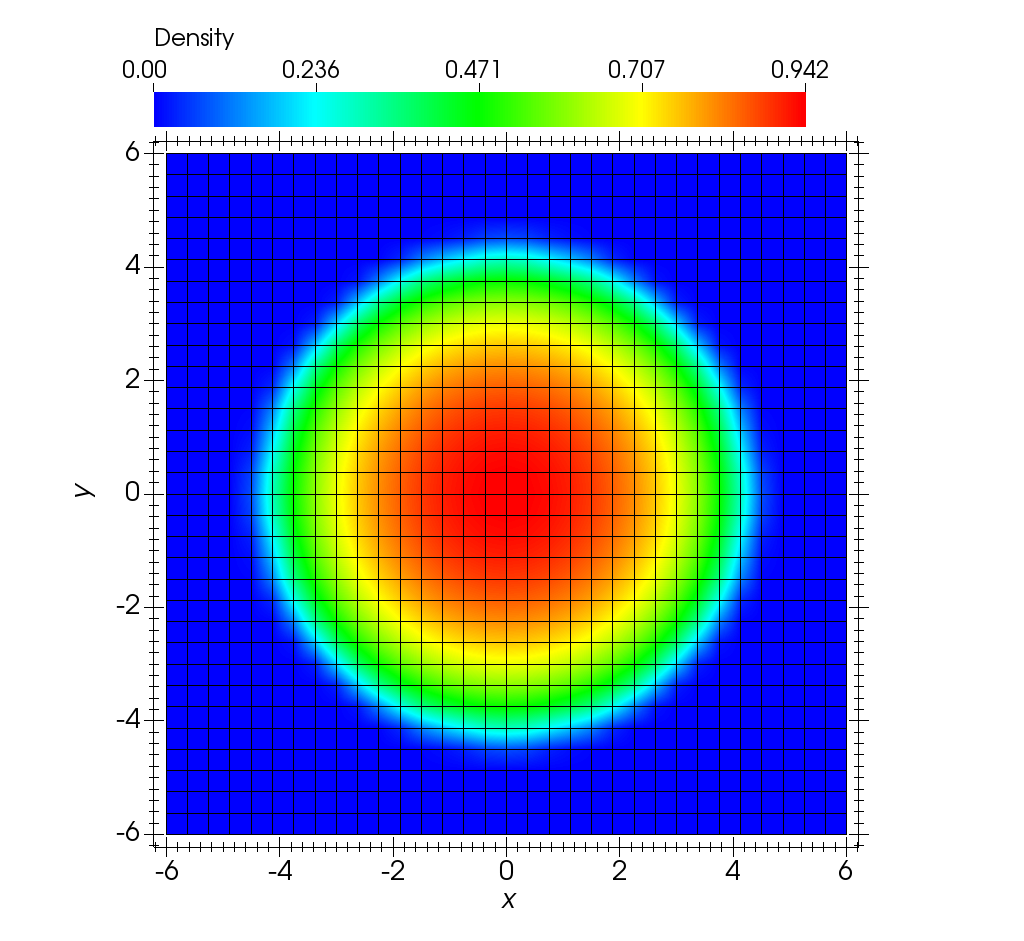}
\captionsetup{justification=centering}\caption{Density; quadrilateral elements}
\end{subfigure}
\begin{subfigure}{.24\textwidth}
\centering
\includegraphics[width=\linewidth,trim={0 0 5.5cm 0},clip]{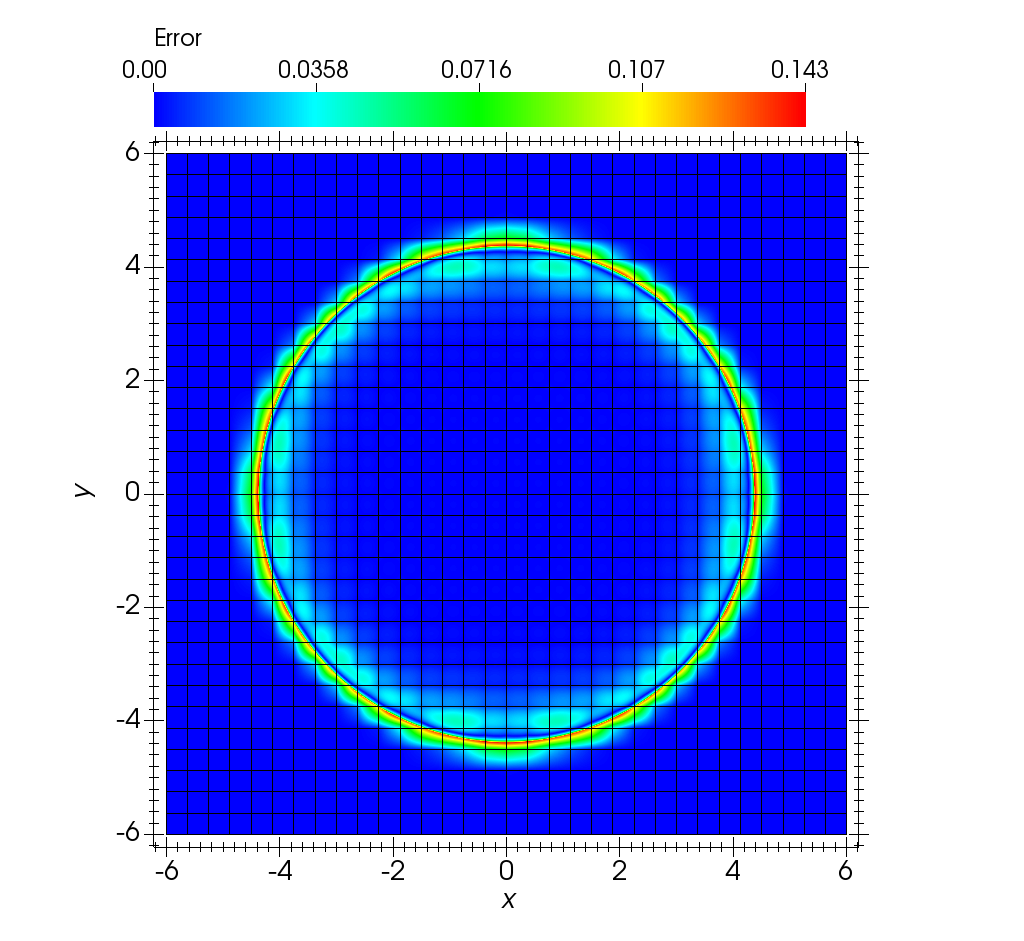}
\captionsetup{justification=centering}\caption{Error; quadrilateral elements}
\end{subfigure}
\small Log-density method
\vspace{2mm}

\begin{subfigure}{.24\textwidth}
\centering
\includegraphics[width=\linewidth,trim={0 0 4cm 0},clip]{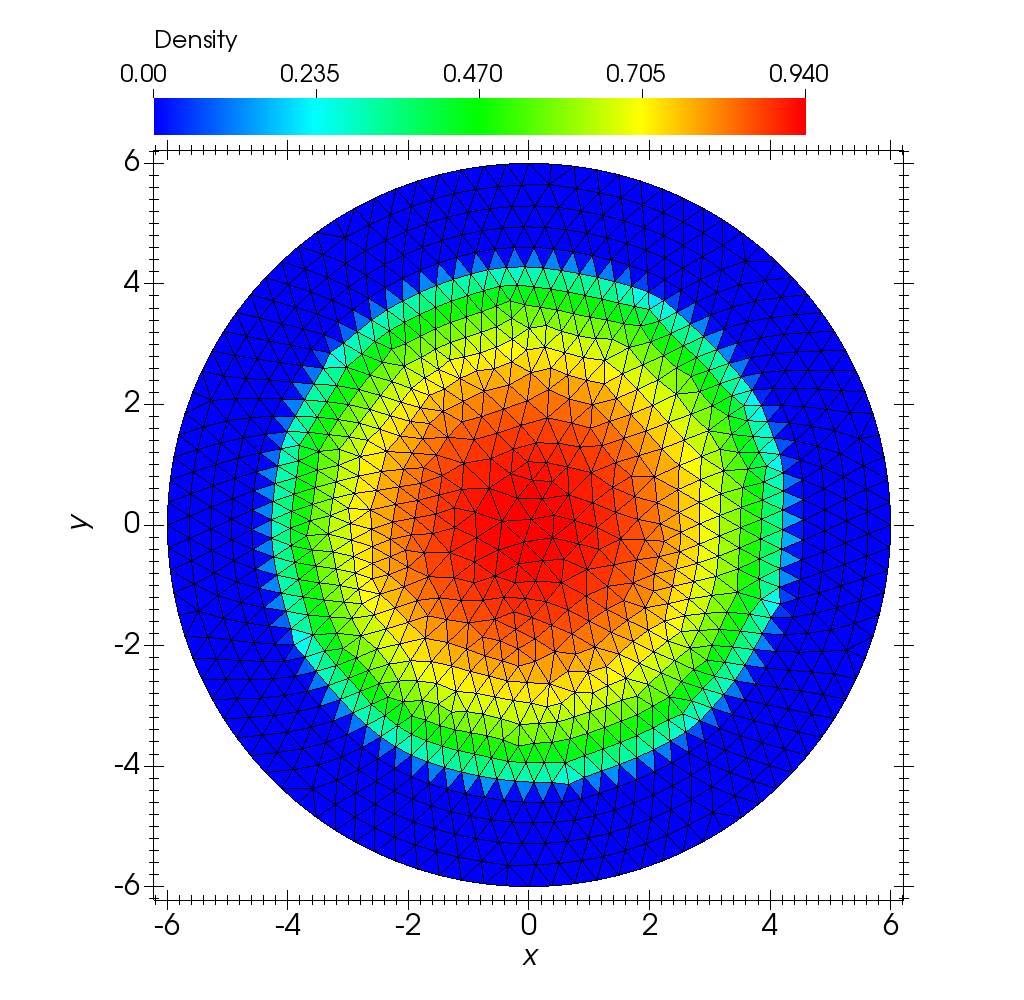}
\captionsetup{justification=centering}\caption{Density; triangular elements}
\end{subfigure}
\begin{subfigure}{.24\textwidth}
\centering
\includegraphics[width=\linewidth,trim={0 0 4cm 0},clip]{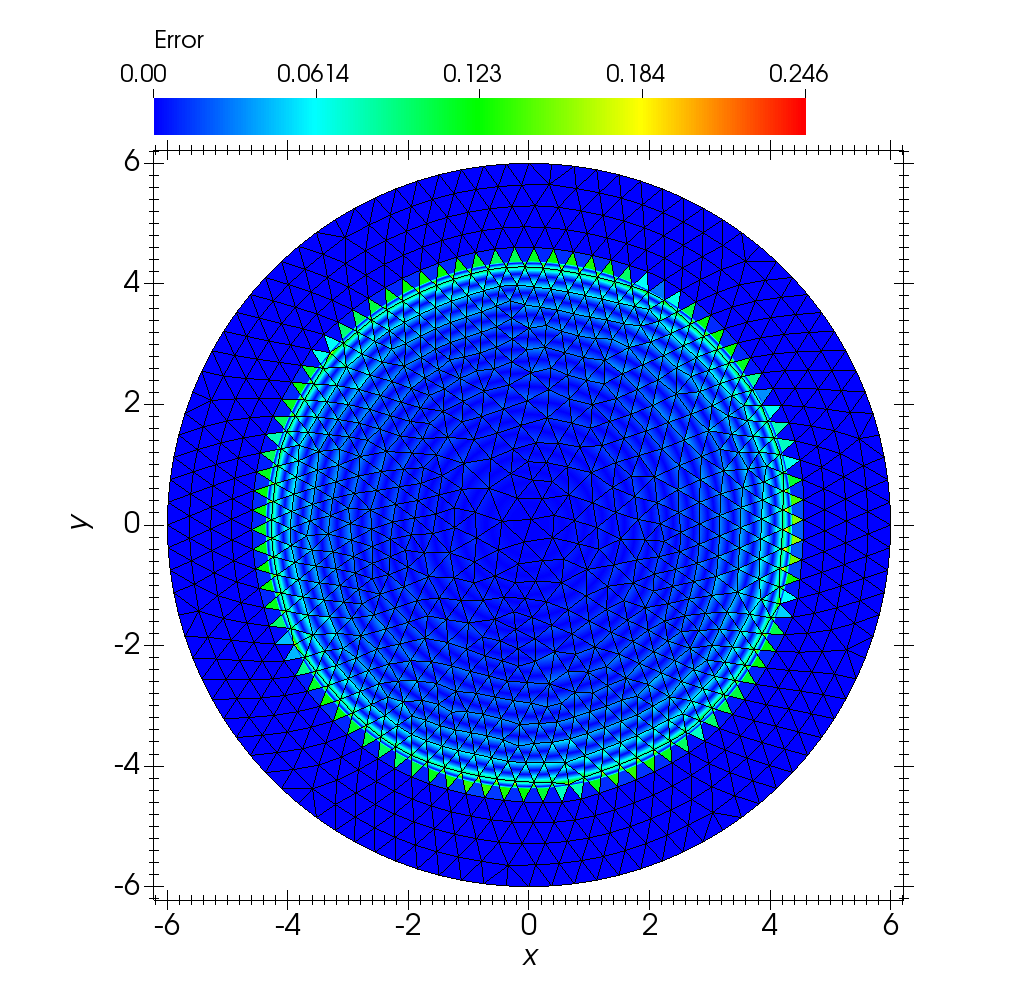}
\captionsetup{justification=centering}\caption{Error; triangular elements}
\end{subfigure}
\begin{subfigure}{.24\textwidth}
\centering
\includegraphics[width=\linewidth,trim={0 0 4cm 0},clip]{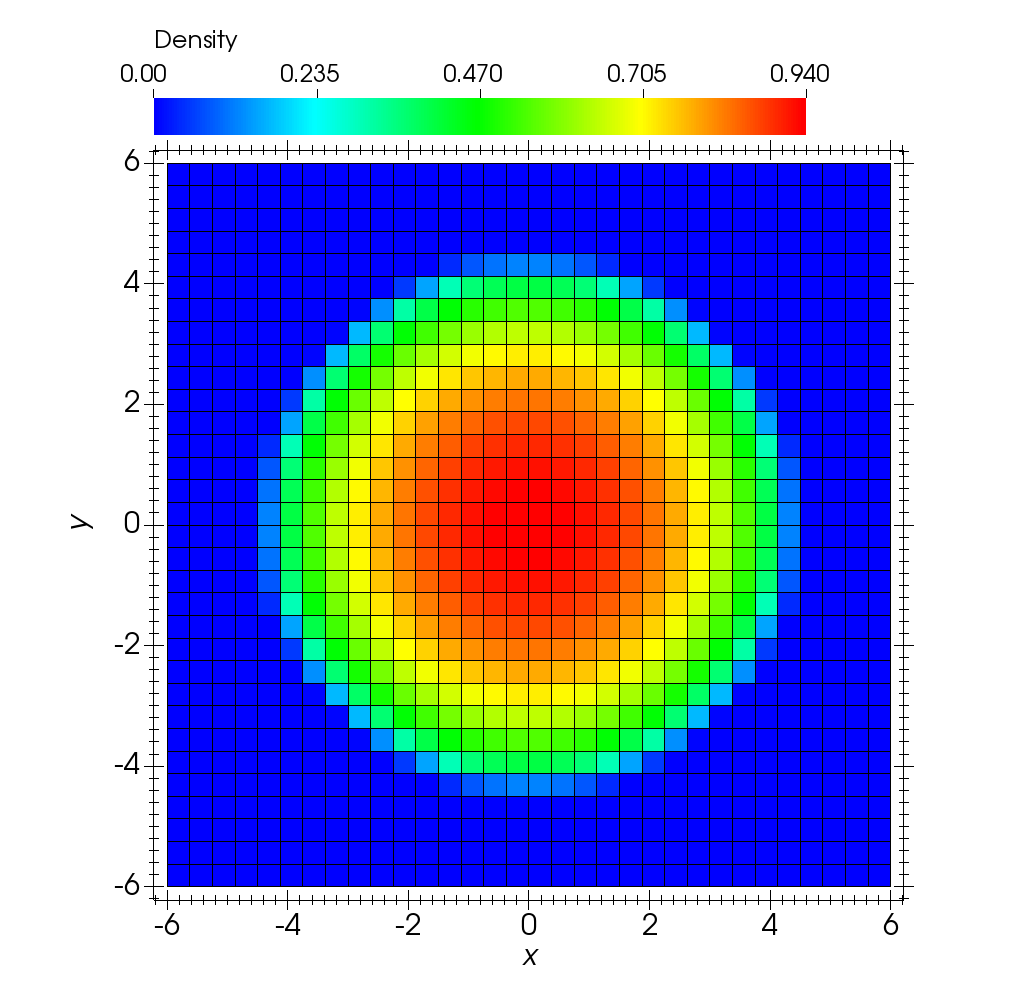}
\captionsetup{justification=centering}\caption{Density; quadrilateral elements}
\end{subfigure}
\begin{subfigure}{.24\textwidth}
\centering
\includegraphics[width=\linewidth,trim={0 0 4cm 0},clip]{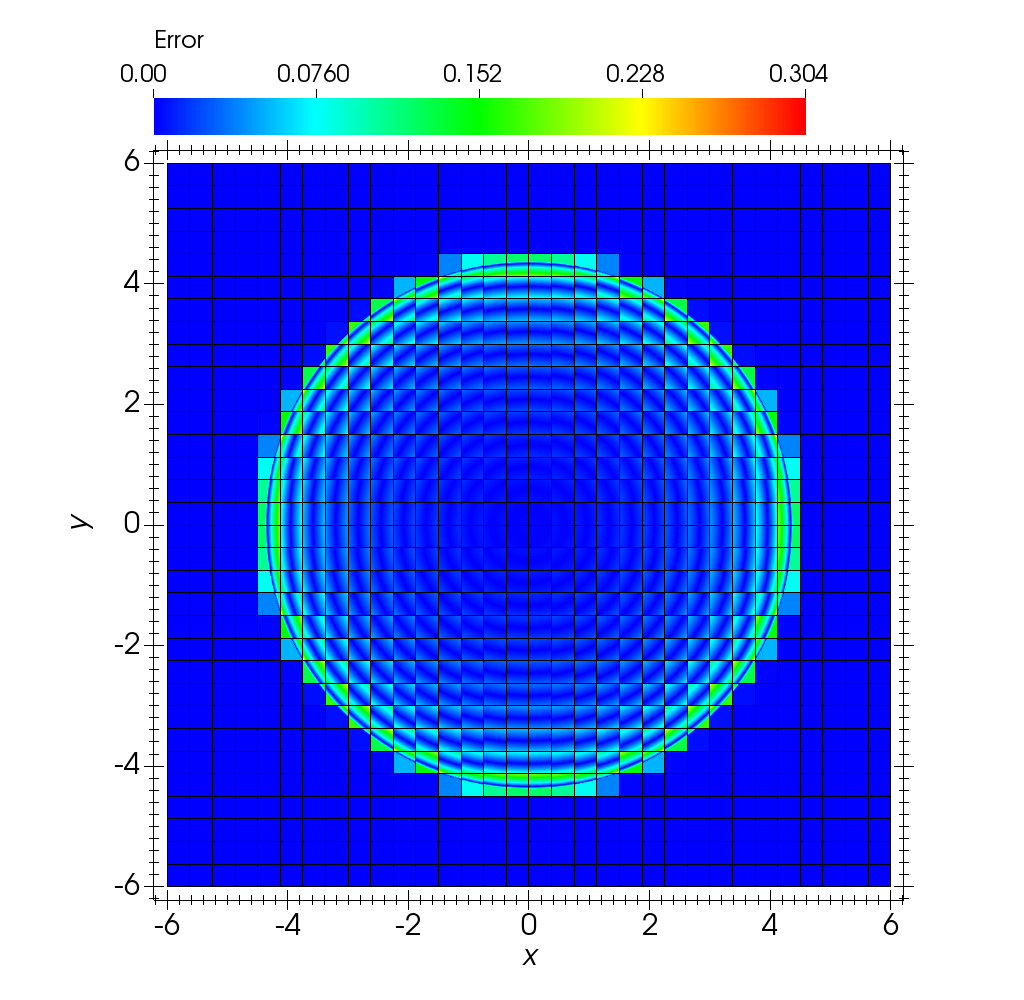}
\captionsetup{justification=centering}\caption{Error; quadrilateral elements}
\end{subfigure}
\small Mixed Method
\vspace{2mm}

\caption{Numerical results at final time $t=1$ for 2D Barenblatt initial data with $m=3$, $\Delta t = 0.025$, and $N\approx 1024$ elements.}
\label{fig:2dbb}
\end{figure} 
\begin{figure}[h!]
\begin{subfigure}{.24\textwidth}
\centering
\includegraphics[width=\linewidth,trim={0 0 5cm 0},clip]{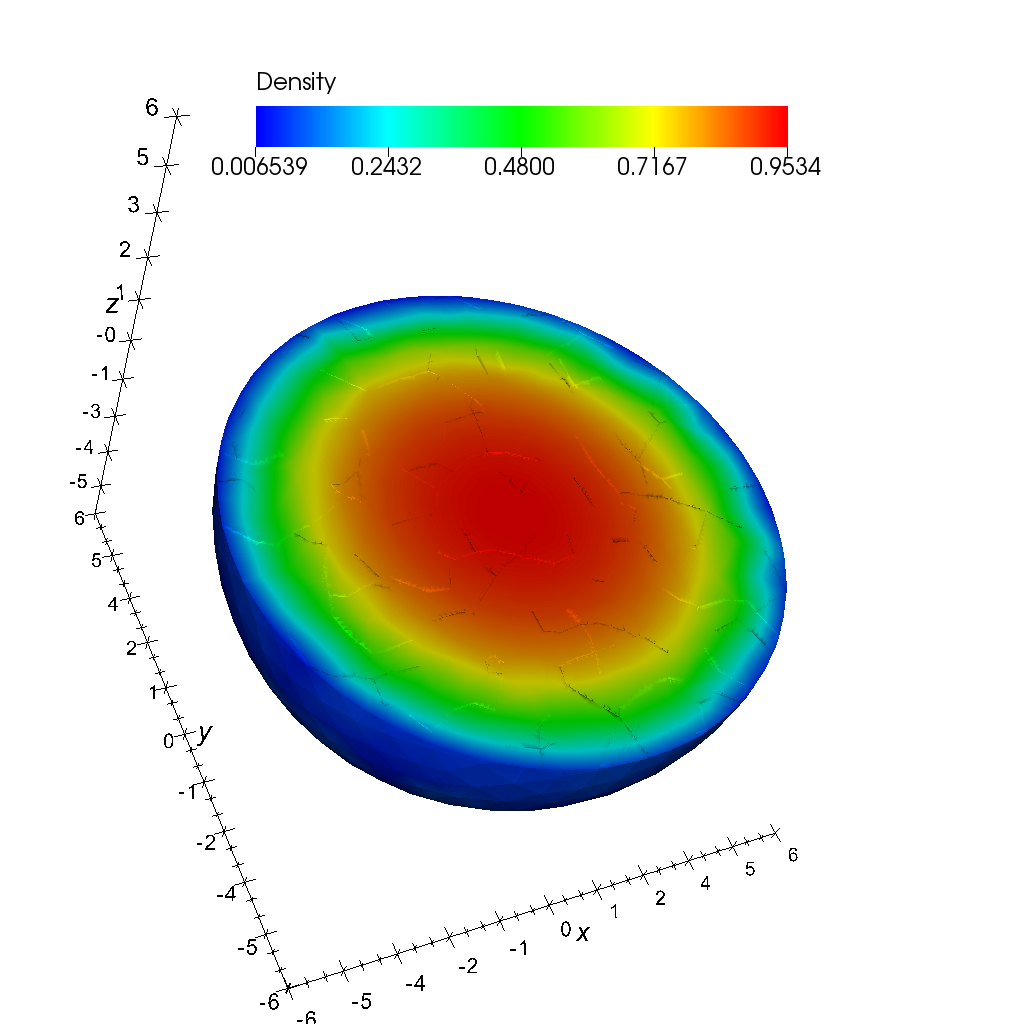}
\captionsetup{justification=centering}\caption{Density; tetrahedral elements.}
\end{subfigure}
\begin{subfigure}{.24\textwidth}
\centering
\includegraphics[width=\linewidth,trim={0 0 5cm 0},clip]{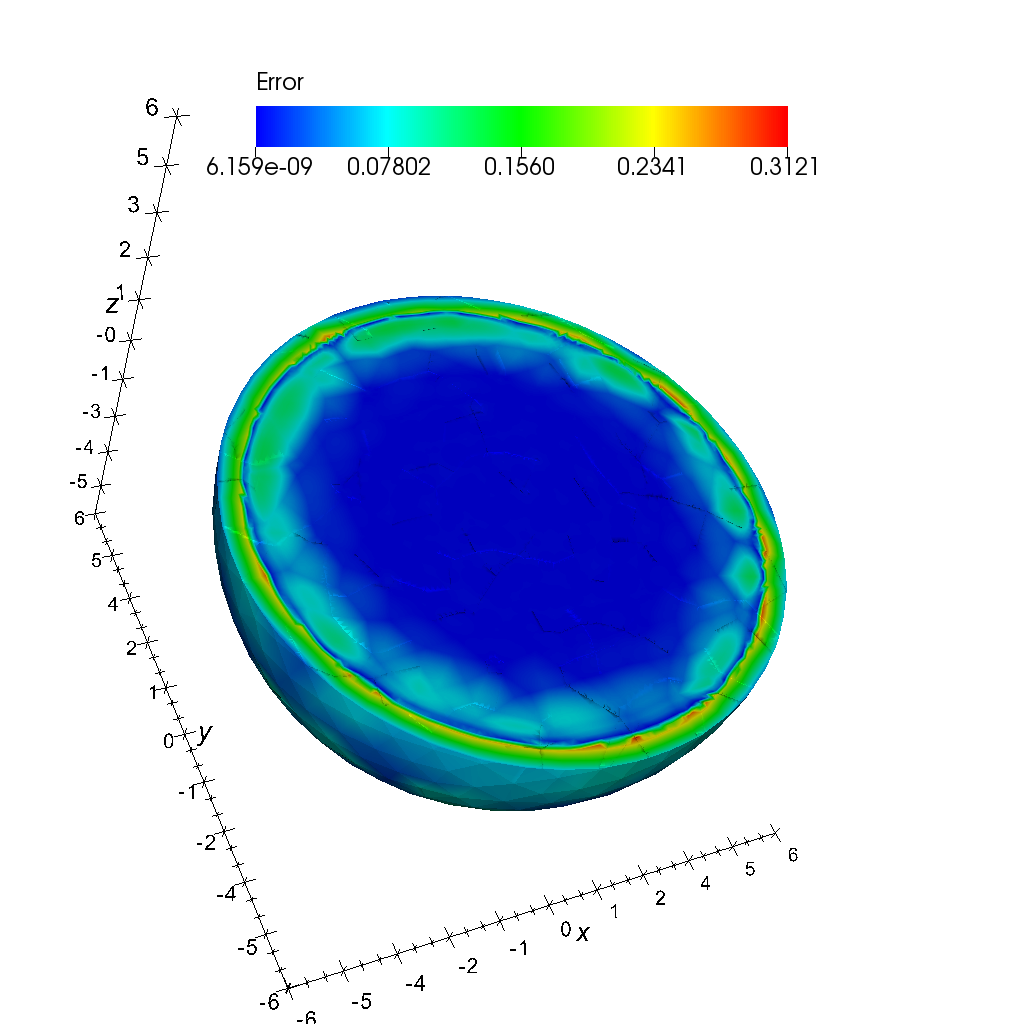}
\captionsetup{justification=centering}\caption{Error; tetrahedral elements}
\end{subfigure}
\begin{subfigure}{.24\textwidth}
\centering
\includegraphics[width=\linewidth,trim={0 0 5cm 0},clip]{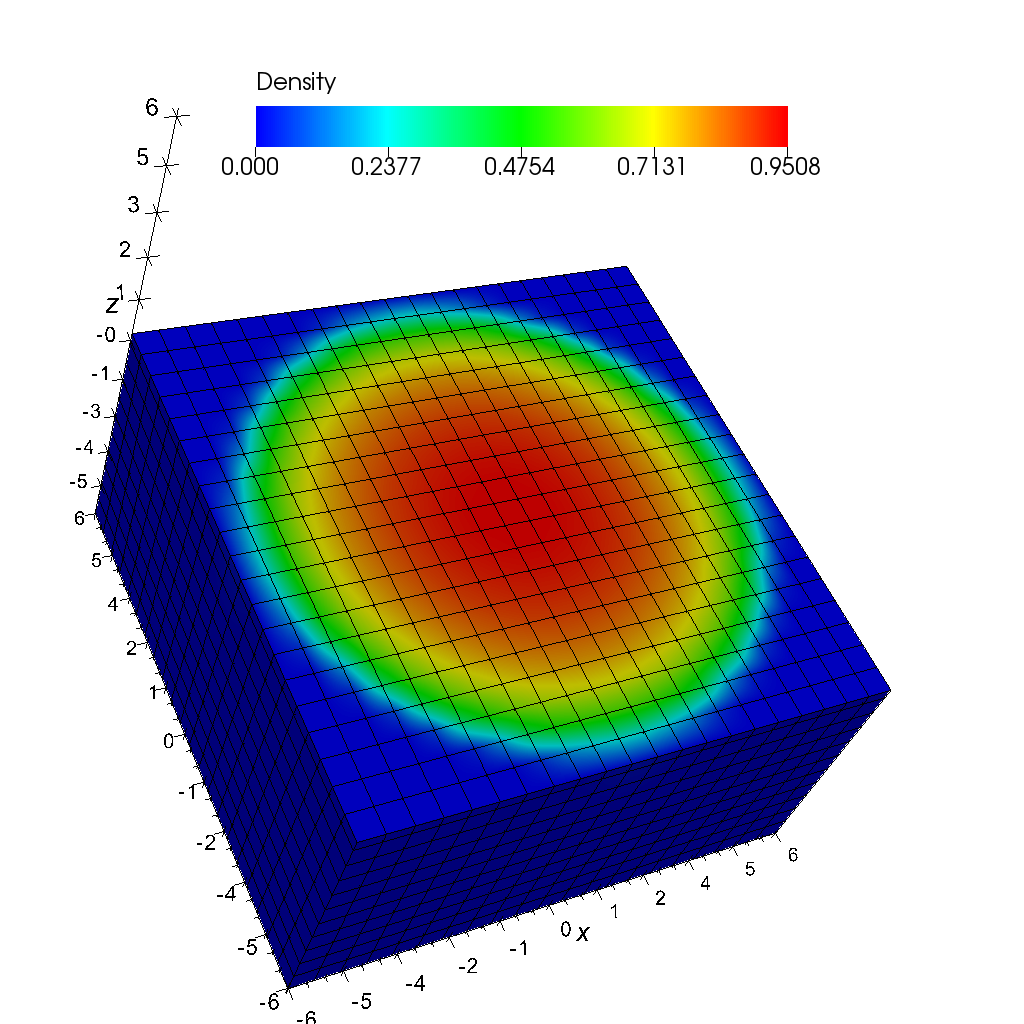}
\captionsetup{justification=centering}\caption{Density; hexahedral elements}
\end{subfigure}
\begin{subfigure}{.24\textwidth}
\centering
\includegraphics[width=\linewidth,trim={0 0 5cm 0},clip]{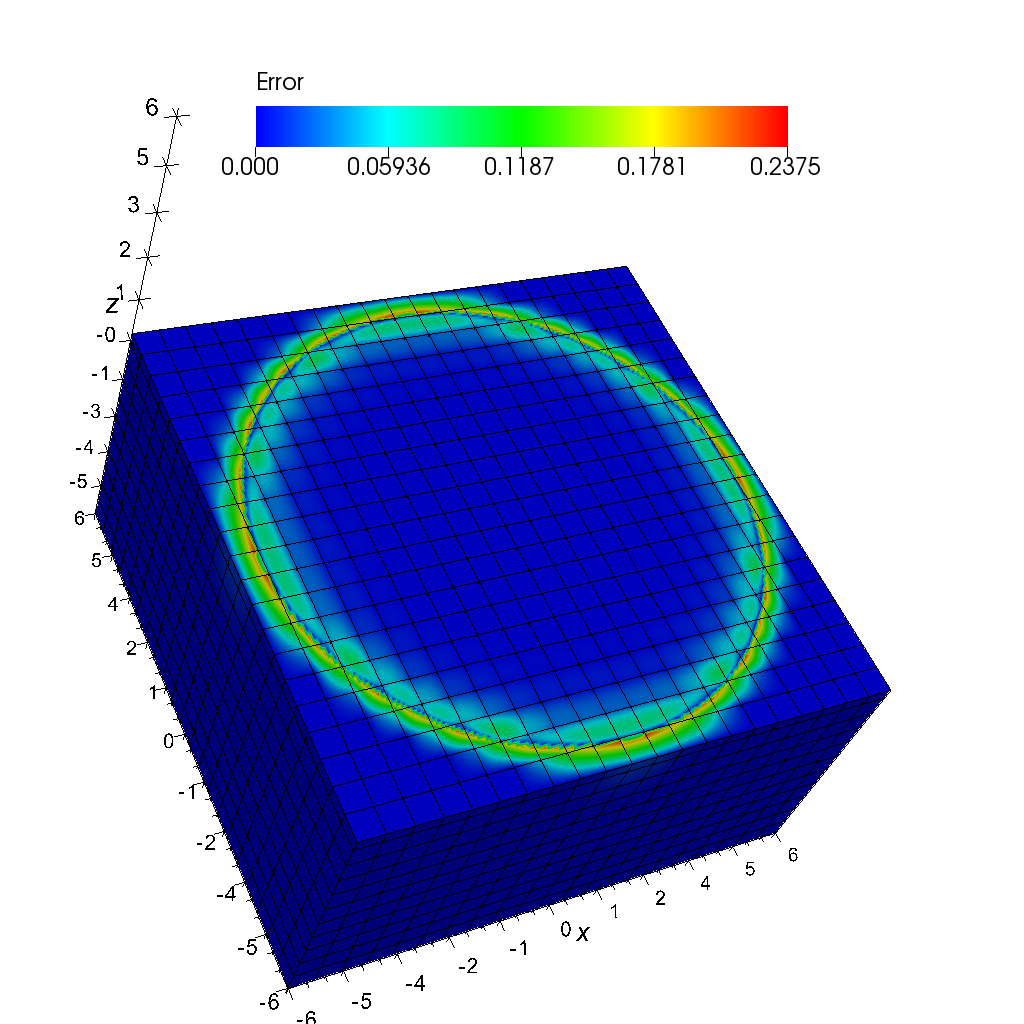}
\captionsetup{justification=centering}\caption{Density; hexahedral elements}
\end{subfigure}
\caption{Numerical results at final time $t=0.2$ using the log-density method for 3D Barenblatt initial data with $m=3$,  $\Delta t = 0.025$, and $N\approx 4100$ elements.}
\label{fig:3dbb}
\end{figure} 
The \textit{d}-dimensional version of \eqref{1dbb} is given by
\begin{equation}
    \label{highdbb}
    \rho_B(\boldsymbol x,t) = (t+1)^{-k} \left(s_0 - \frac{k(m-1)}{2dm} \frac{|\boldsymbol x|^2}{(t+1)^{2k/d}}\right)_+^{\frac{1}{m-1}}, \quad t > 0.
\end{equation}
Analogous to the one dimensional case, we verify the accuracy of the two schemes in two dimensions by using \eqref{highdbb} as initial data with $x \in [-6,6]^2$, $s_0=1$, and $m\in \{2,3,4\}$ and conduct a spacetime mesh-refinement study for the $L^2$ error at final time $t=0.2$. The $L^2$-error is recorded in the entire domain and in the box $[-3,3]^2$ away from the interface. We use a sequence of spatial meshes with $\{(32\times 32)\cdot 2^i\}_{i=0}^3$ tensor-product elements and set the timestep size to $\frac{1}{5\cdot 4^{i}}$ for the log-density method and to $\frac{1}{10\cdot 2^{i}}$ for the mixed method, correspondingly. 
The results of the study are summarized in Table \ref{t:2DBB}. Akin to the one dimensional case, we observe second-order spatial accuracy and first-order temporal accuracy of the log-density method in the box $[-3,3]^2$ where the solution remains smooth at final time. In the same region, the mixed method is observed to be first-order accurate in space and time as expected. Additionally, the order of convergence is observed to decay with $m$ when the error is measured in the entire domain. This deterioration is slower for the case of the mixed method.

In figure \ref{fig:2dbb}, we plot the density and error profiles obtained using the two schemes at final time $t=0.2$ for the data \eqref{highdbb} with $m=3$, $\Delta t = 0.025$, and $N\approx 1024$ elements. Results on both triangular and quadrilateral meshes are shown. Panels (A)-(D) display the profiles computed using the log-density method, while (E)-(F) display the profiles computed using the mixed method. 

We additionally evolve the three-dimensional version of \ref{highdbb} on tetrahedral and hexahedral meshes of the domain $\Omega_h = [-6,6]^3$ using the log-density method until final time $t=0.2$ with $\Delta t=0.025$ and $N \approx 4100$ elements. We set $m=3$ and $s_0=1$. The plots of this simulation are shown in figure \ref{fig:3dbb}. In the region $[-3,3]^3$, the $L^2$-error is 0.0196 for hexahedral elements and 0.0674 for tetrahedral elements. We remark that the mixed method can also be used for this case, but only hexahedral or regular tetrahedral meshes may be used due to the lack of a proper mass-lumping quadrature formula for general tetrahedral elements.

\subsection{Merging Gaussians}
\begin{figure}[]
\begin{subfigure}{.3\textwidth}
\centering
\includegraphics[width=\linewidth,trim={0 0 3cm 0},clip]{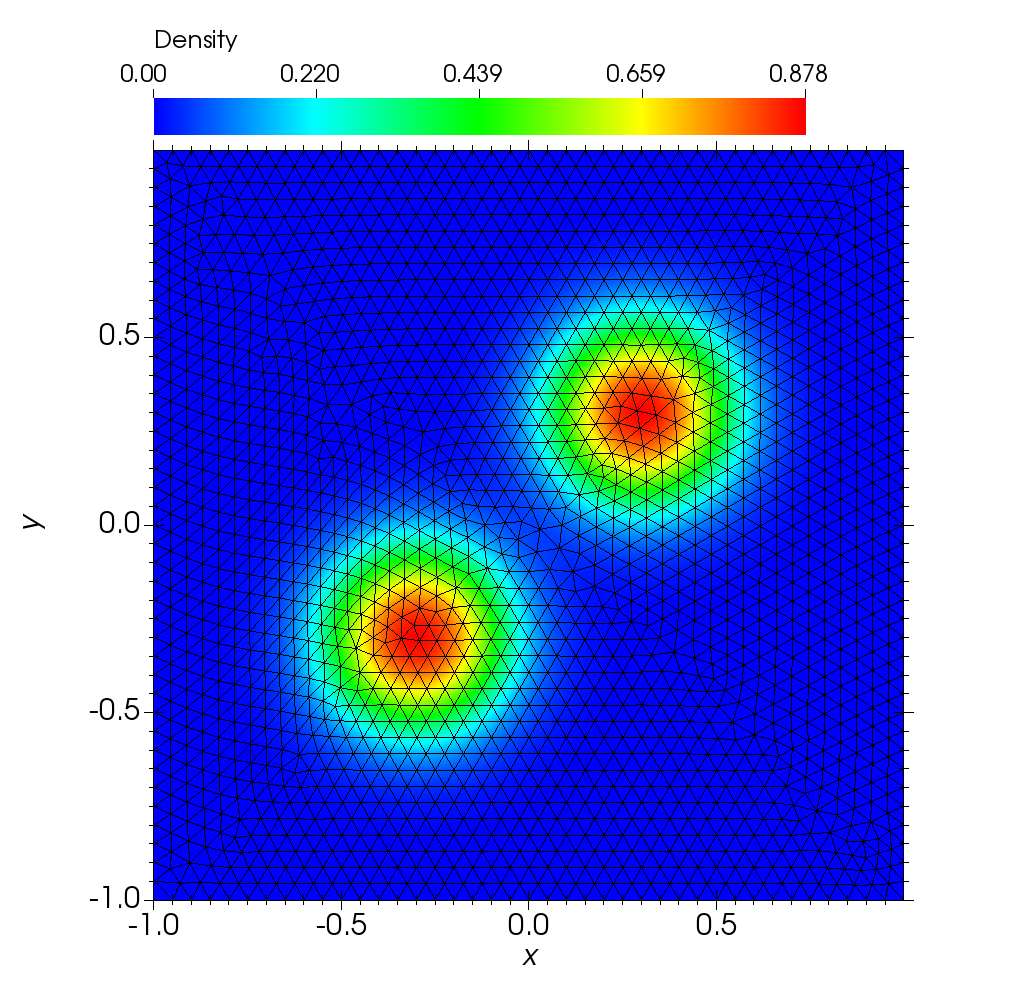}
\caption{$t=0$}
\end{subfigure}
\begin{subfigure}{.3\textwidth}
\centering
\includegraphics[width=\linewidth,trim={0 0 3cm 0},clip]{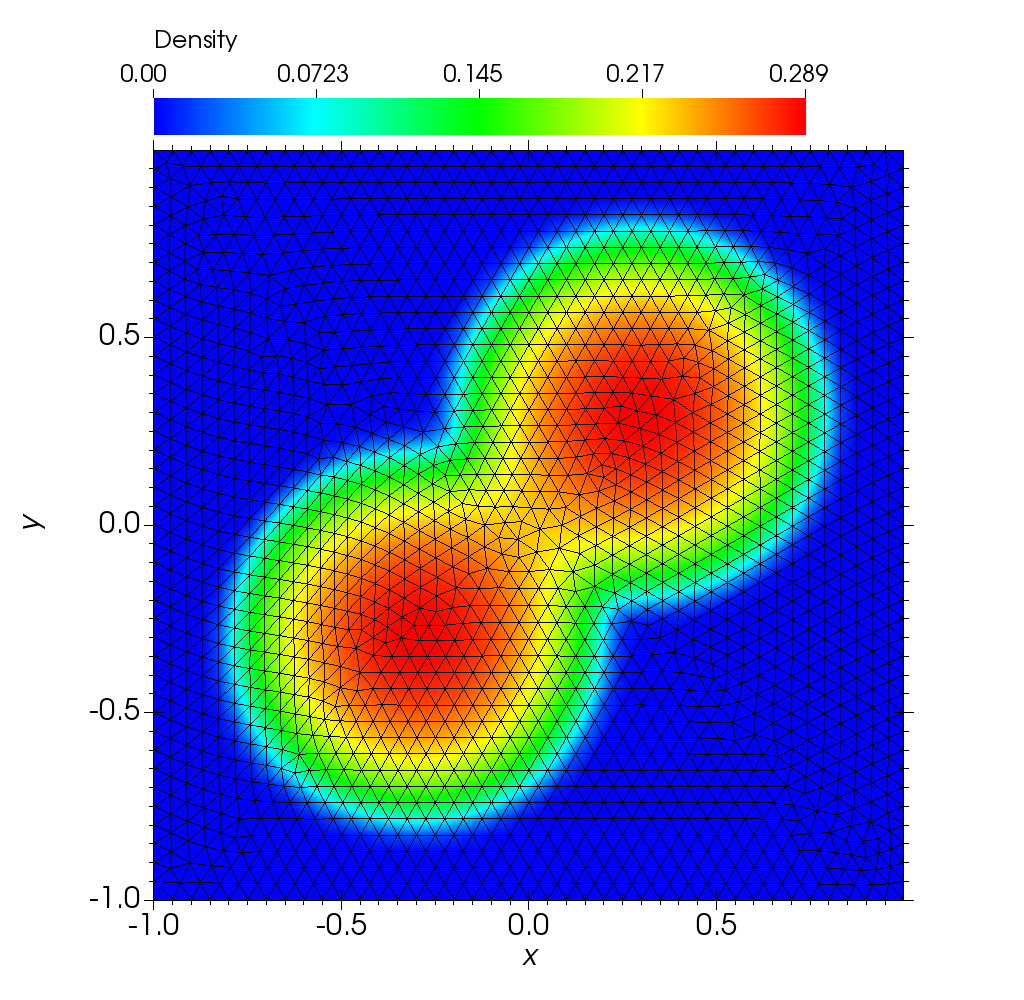}
\caption{$t=0.15$}
\end{subfigure}
\begin{subfigure}{.3\textwidth}
\centering
\includegraphics[width=\linewidth,trim={0 0 3cm 0},clip]{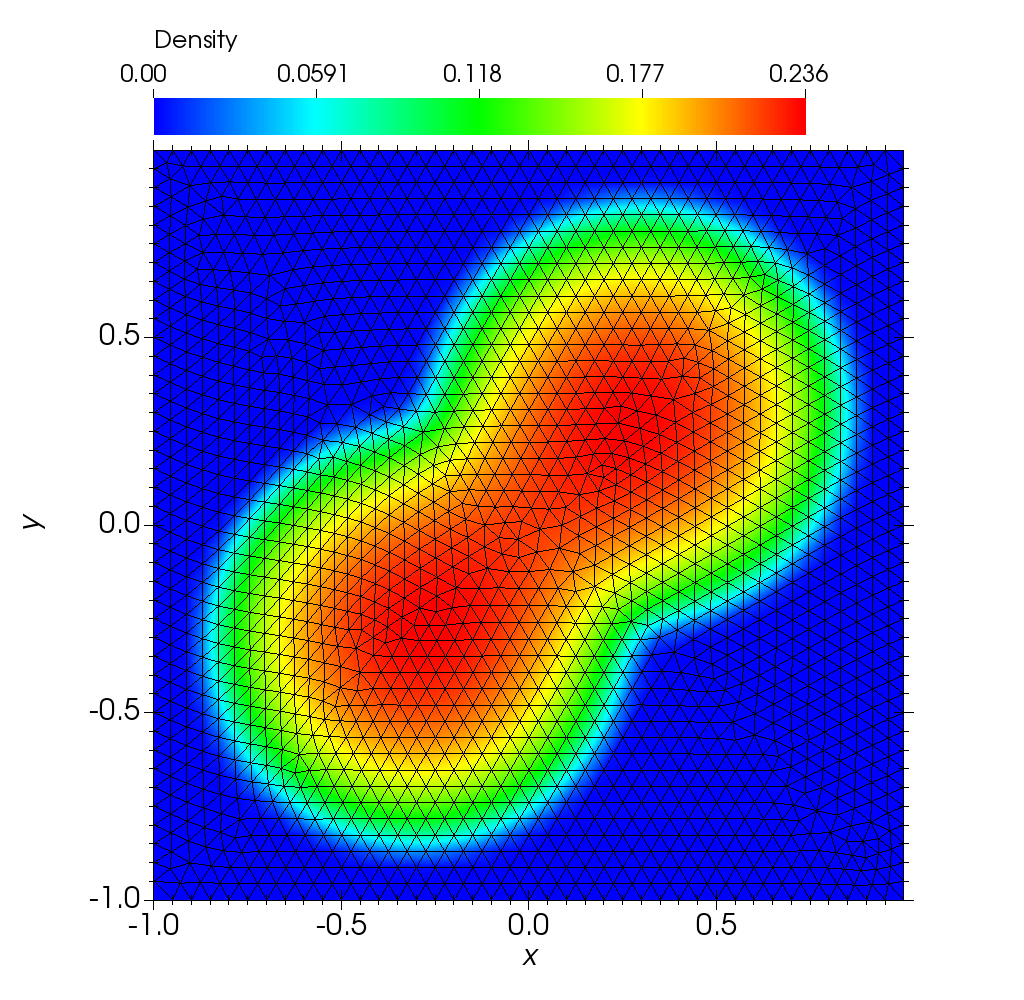}
\caption{$t=0.3$}
\end{subfigure}
\small Log-density method
\vspace{2mm}

\begin{subfigure}{.3\textwidth}
\centering
\includegraphics[width=\linewidth,trim={0 0 3cm 0},clip]{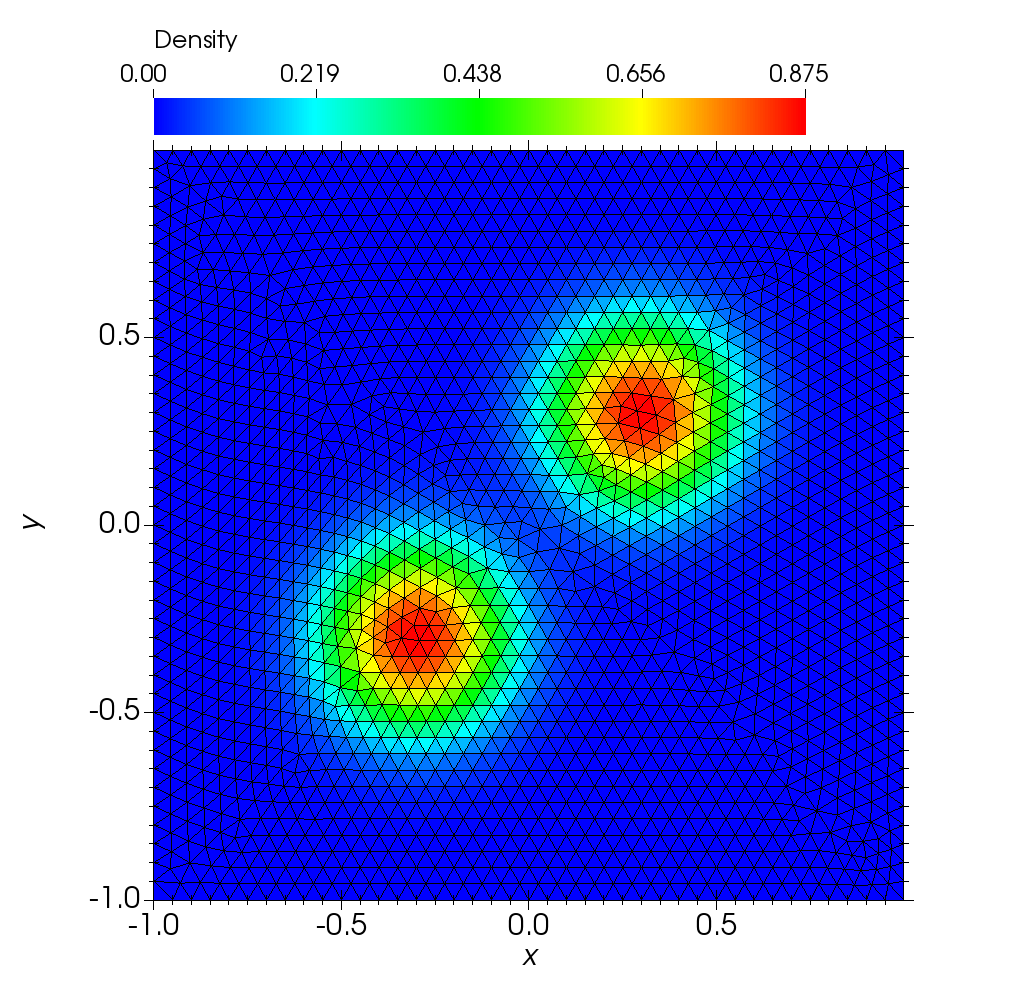}
\caption{$t=0$}
\end{subfigure}
\begin{subfigure}{.3\textwidth}
\centering
\includegraphics[width=\linewidth,trim={0 0 3cm 0},clip]{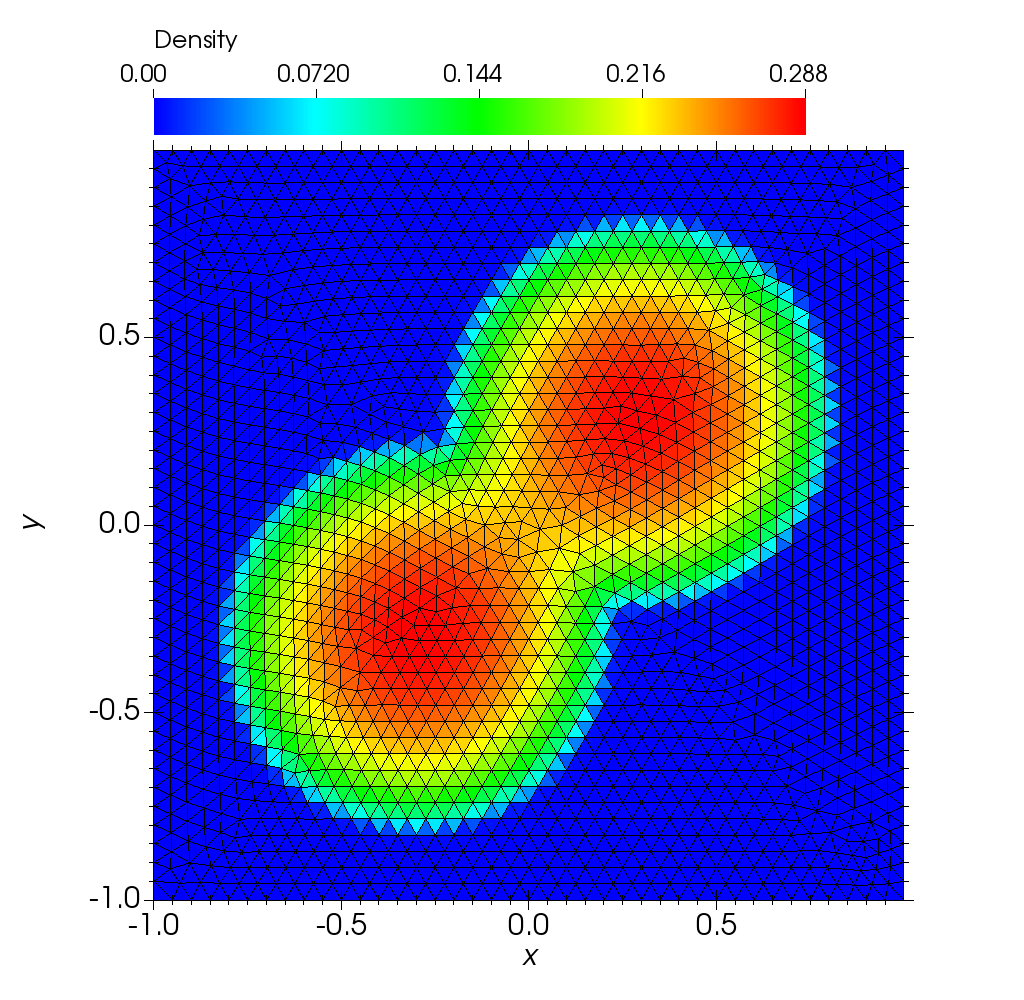}
\caption{$t=0.15$}
\end{subfigure}
\begin{subfigure}{.3\textwidth}
\centering
\includegraphics[width=\linewidth,trim={0 0 3cm 0},clip]{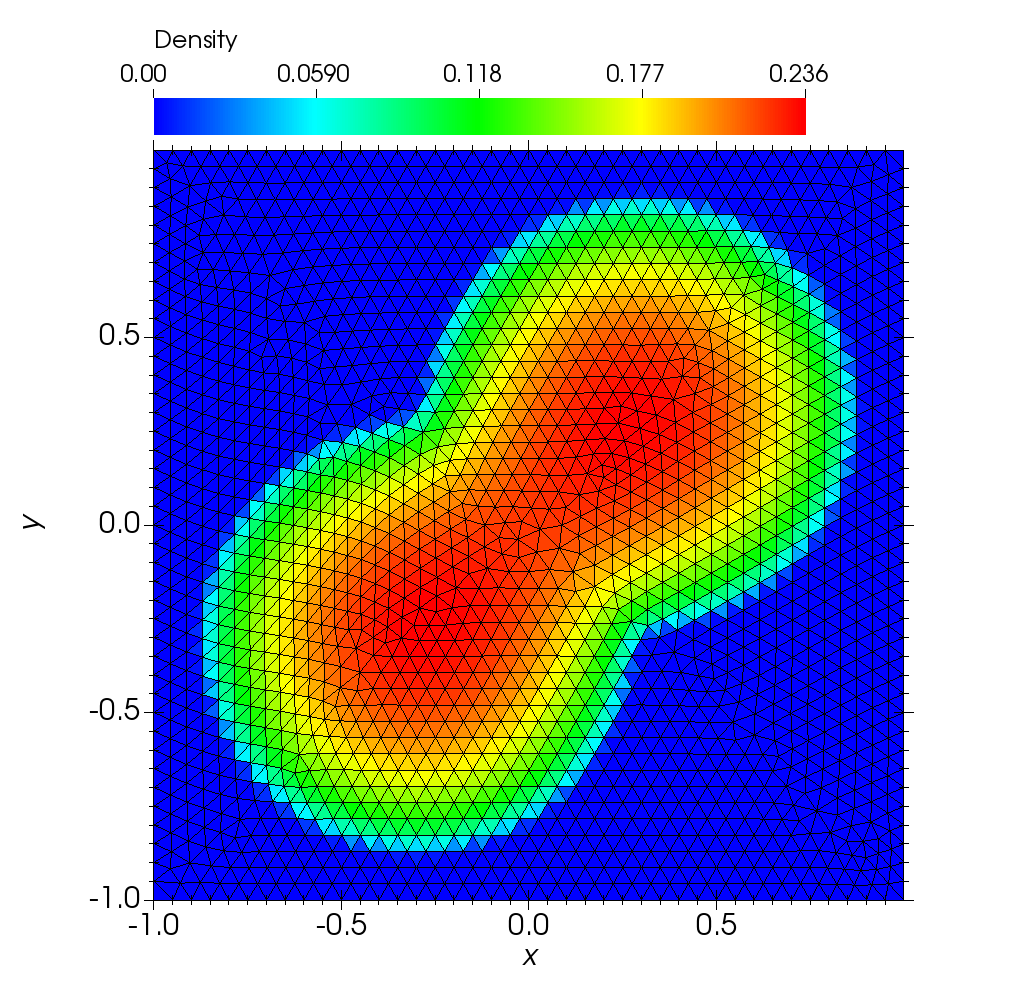}
\caption{$t=0.3$}
\end{subfigure}

\small Mixed method
\vspace{2mm}
\caption{Merging Gaussians test with $ m =3, \Delta t = 0.001$, and $N = 3750$ triangular elements.}
\label{fig:mg}
\end{figure}

We further investigate the robustness of the two schemes on a popular test case consisting of two initial Gaussian peaks that merge into a single peak under the action of the PME \cite{Liu20a,Ngo17,Carrillo}. The initial condition for this test is given by
\begin{equation}
\label{mg}
    \rho^0(x, y)=e^{-20\left((x-0.3)^2+(y-0.3)^2\right)}+e^{-20\left((x+0.3)^2+(y+0.3)^2\right)},
\end{equation}
\noindent where the domain is taken to be $\Omega = [-1,1]^2$. We set $m=3$, $\Delta t = 0.001$, and evolve the initial data \eqref{mg} using the two schemes on a triangular grid with 3750 elements. We record the solution at $t=0$, $t=0.15$, and $t=0.3$, and plot the results of the simulation in figure \ref{fig:mg}, where one can observe that in both cases the two peaks move towards each other and eventually start merging.

\subsection{Complex support}
\begin{figure}[]
\begin{subfigure}{.3\textwidth}
\centering
\includegraphics[width=\linewidth]{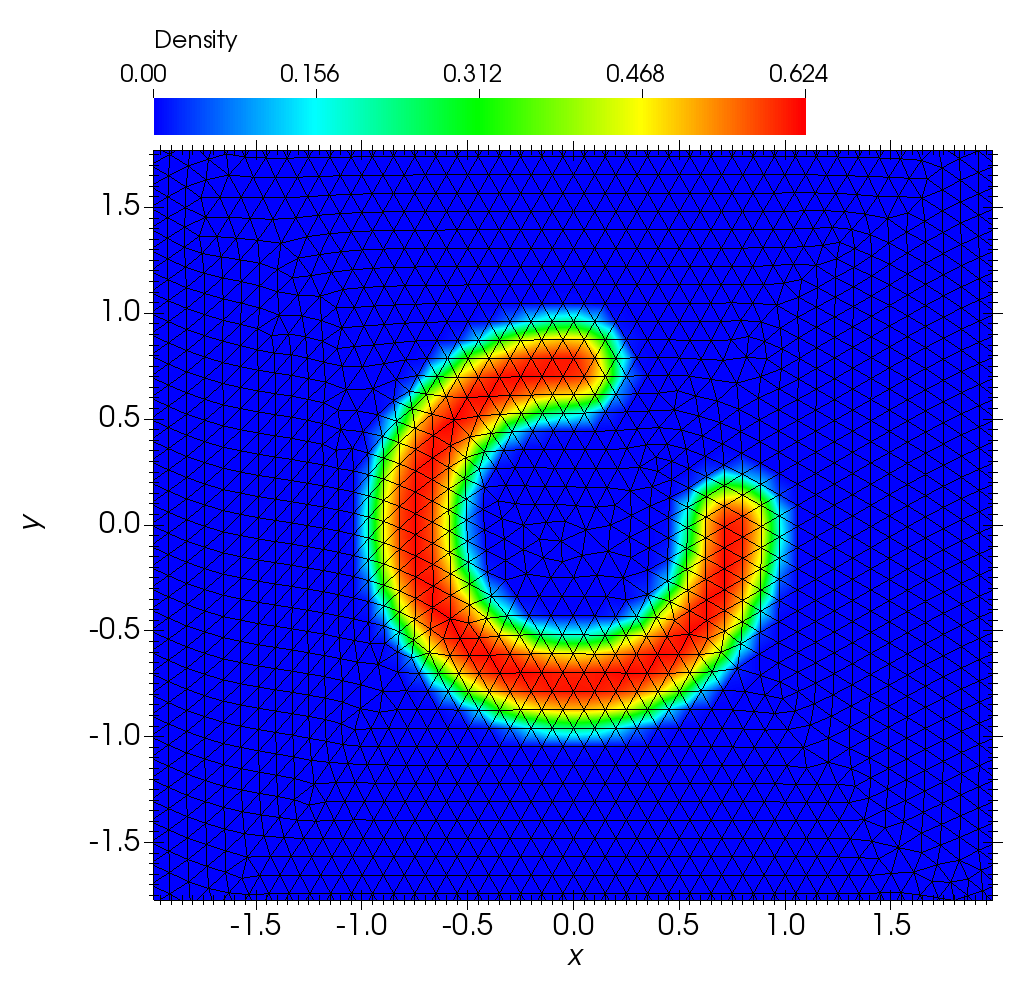}
\caption{$t=0$}
\end{subfigure}
\begin{subfigure}{.3\textwidth}
\centering
\includegraphics[width=\linewidth]{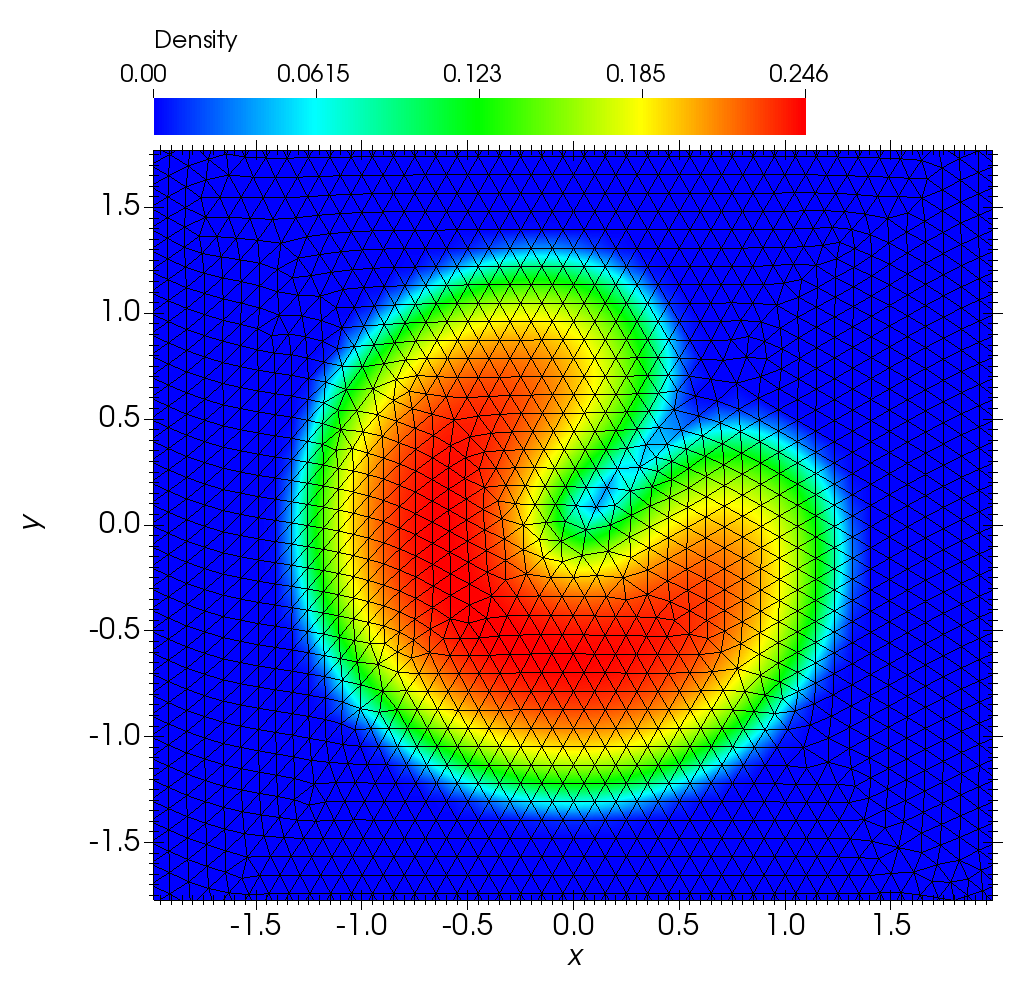}
\caption{$t=0.5$}
\end{subfigure}
\begin{subfigure}{.3\textwidth}
\centering
\includegraphics[width=\linewidth]{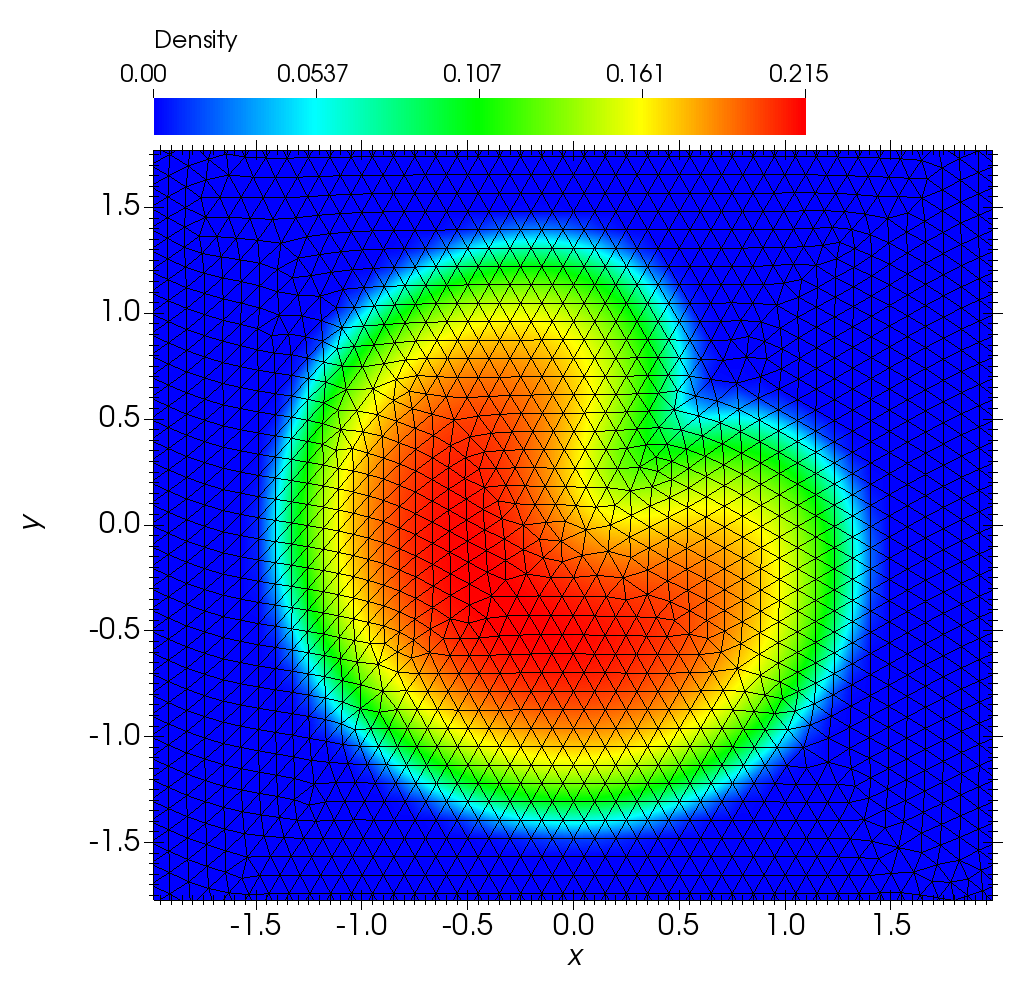}
\caption{$t=1.0$}
\end{subfigure}
\small Log-density method
\vspace{2mm}

\begin{subfigure}{.3\textwidth}
\centering
\includegraphics[width=\linewidth]{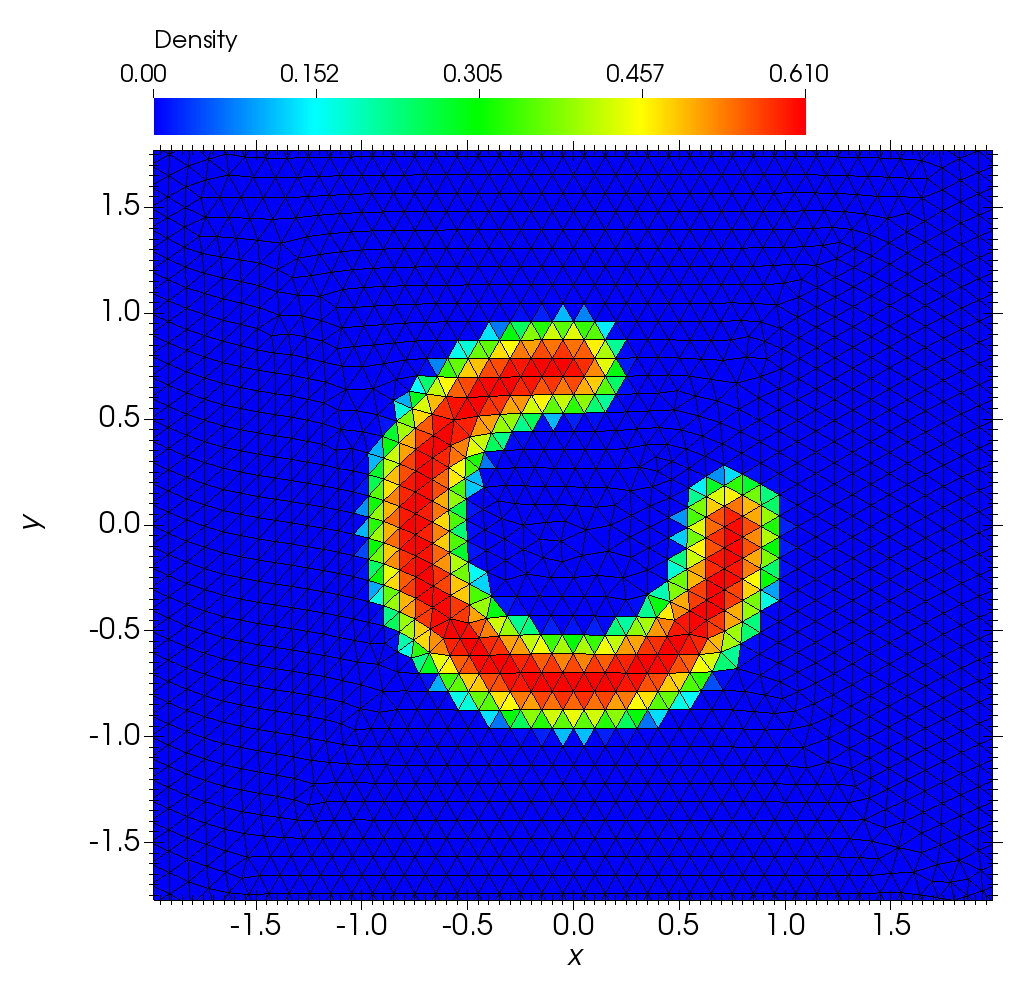}
\caption{$t=0$}
\end{subfigure}
\begin{subfigure}{.3\textwidth}
\centering
\includegraphics[width=\linewidth]{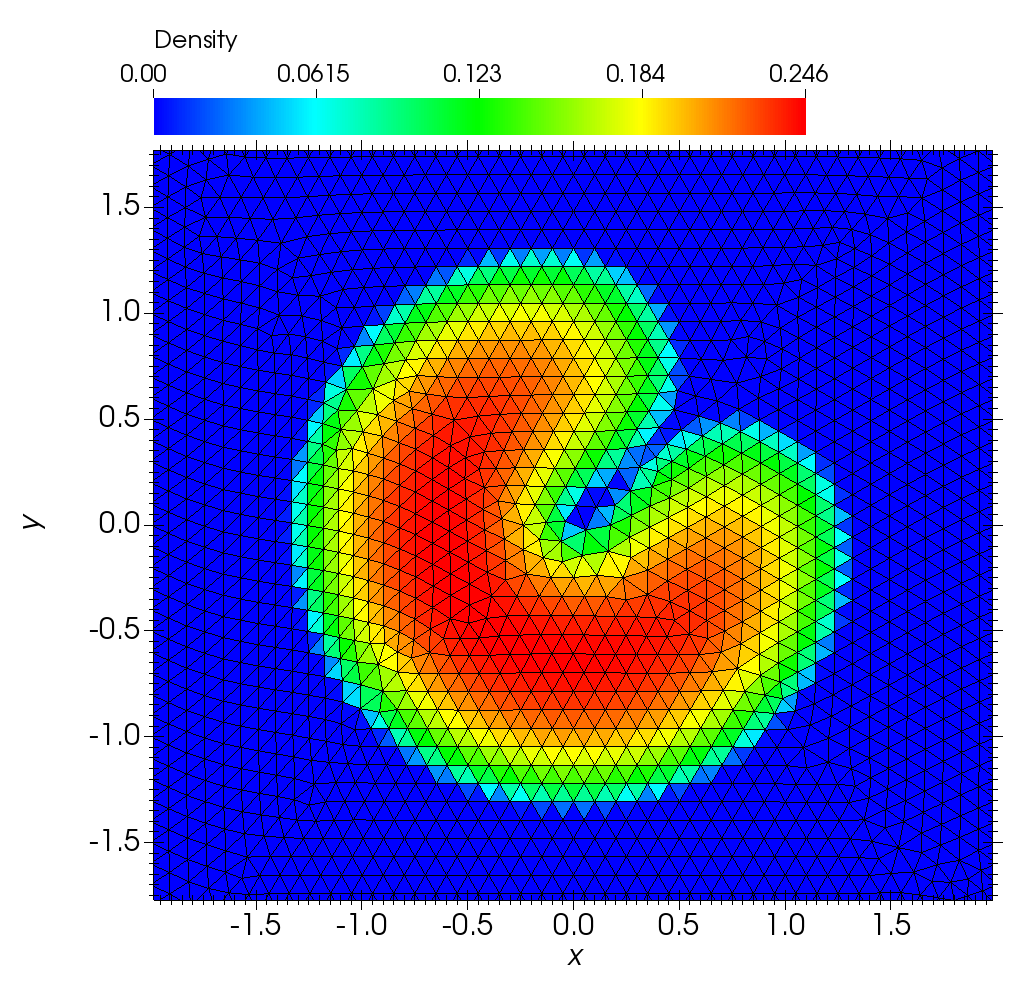}
\caption{$t=0.5$}
\end{subfigure}
\begin{subfigure}{.3\textwidth}
\centering
\includegraphics[width=\linewidth]{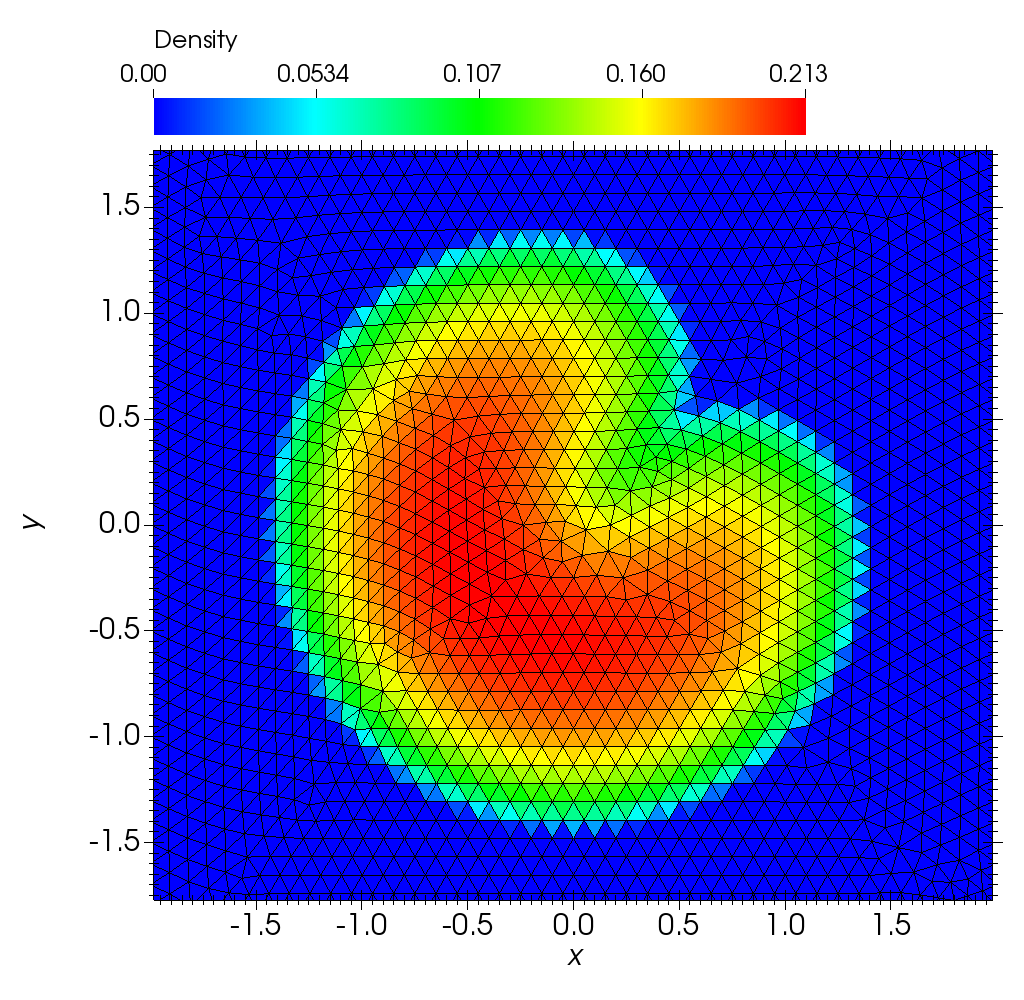}
\caption{$t=1.0$}
\end{subfigure}

\small Mixed method
\vspace{2mm}
\caption{Complex support test with $ m =3, \Delta t = 0.001$, and $N= 3750$ triangular elements.}
\label{fig:cs}
\end{figure}
In the final numerical example, we consider the following initial data \cite{Liu20a,ngo2019,BAINES2005450}
\small
\begin{equation}
\label{cs}
    \rho^0(x, y)=
\left\{\begin{aligned}
25\left(0.25^2-\left(\sqrt{x^2+y^2}-0.75\right)^2\right)^{\frac{3}{2(m-1)}}, &\quad \sqrt{x^2+y^2} \in[0.5,1] \text { and }(x<0 \text { or } y<0), \\
25\left(0.25^2-x^2-(y-0.75)^2\right)^{\frac{3}{2(m-1)}}, &\quad x^2+(y-0.75)^2 \leq 0.25^2 \text { and } x \geq 0, \\
25\left(0.25^2-(x-0.75)^2-y^2\right)^{\frac{3}{2(m-1)}}, &\quad(x-0.75)^2+y^2 \leq 0.25^2 \text { and } y \geq 0, \\
0, &\quad \text { otherwise},
\end{aligned}\right.
\end{equation}
\normalsize
on the domain $[-2,2]$. The support of intial data \eqref{cs} has the shape of a horseshoe or a partial donut. We set $m=3$, $\Delta t = 0.001$, and evolve the initial data using the two schemes on a triangular grid consisting of 3750 elements. The solution is recorded at times $t=0$, $t=0.5$, and $t=1.0$, and the computed profiles are presented in figure \ref{fig:cs}. The plots demonstrate that both the schemes deliver reliable results, as the horseshoe ends are observed to evolve towards each other before ultimately intersecting. Additionally, in line with previous observations, the boundary of the support is captured more sharply by the mixed method. 

We note that unlike Lagarangian schemes, the two methods exhibit a robust handling of the topology change without requiring any interpolation of the solution on to a new mesh when this event occurs. 

\section{Conclusion}
We presented two distinct, first-order spacetime accurate, finite element approaches to the PME. The log-density approach is constructed for a problem with Neumann boundary conditions, and the properties of mass conservation, energy stability, unique solvability, and bound preservation on unstructured Delaunay meshes are proved. The scheme is shown to be second order in space and first order in time. The mixed approach is constructed for a problem with Dirichlet boundary conditions and is shown to be mass conservative, energy stable, and positivity preserving (under a CFL condition). The mixed scheme is shown to be first order in both space and time. Both schemes can handle compactly supported initial data without the need of any perturbation.

\

\noindent\textbf{Funding} This research is partially supported by NSF grant DMS-2012031.

\

\noindent\textbf{Data Availability} Enquiries about data availability shall be directed to the authors.

\section*{Declaration}
\noindent\textbf{Conflict of interest}
The authors declare that there are no known conflicts of interest assocaited with this work.

\bibliography{fsi}
\bibliographystyle{siam}

\end{document}